\pdfoutput=1

\documentclass[11pt]{article}

\usepackage[T1]{fontenc}
\usepackage[utf8]{inputenc}
\usepackage{lmodern}
\usepackage{microtype}
\usepackage[a4paper,hmargin=20mm,vmargin=25mm]{geometry}
\usepackage[fleqn]{amsmath}
\usepackage{amssymb}
\usepackage{graphicx}
\usepackage{xcolor}
\usepackage{placeins}
\usepackage{booktabs}
\usepackage{tabularx}
\usepackage{array}
\usepackage{soul}
\usepackage{makecell}
\usepackage{multirow}
\usepackage[round,authoryear,longnamesfirst]{natbib}
\usepackage{authblk}
\usepackage{fancyhdr}
\usepackage[colorlinks=true,linkcolor=blue!60!black,citecolor=blue!60!black,urlcolor=blue!60!black]{hyperref}
\hypersetup{
  pdftitle={The r-Safety Reserve and Adaptive Kinematic Smoothing for Car-Following Transitions},
  pdfauthor={Xuesong (Simon) Zhou, Ziyi Zhang}}

\title{\textbf{The $r$-Safety Reserve and Adaptive Kinematic Smoothing for Car-Following Transitions}}

\author[1]{Xuesong (Simon) Zhou\thanks{\href{mailto:xzhou74@asu.edu}{xzhou74@asu.edu}; ORCID: \href{https://orcid.org/0000-0002-9963-5369}{0000-0002-9963-5369}}}
\author[1]{Ziyi Zhang\thanks{\href{mailto:zzhan657@asu.edu}{zzhan657@asu.edu}; ORCID: \href{https://orcid.org/0000-0001-8720-9840}{0000-0001-8720-9840}}}
\affil[1]{School of Sustainable Engineering and the Built Environment, Arizona State University, Tempe, AZ 85281, USA}

\date{}

\begin{document}

\maketitle
\thispagestyle{empty}

\begin{abstract}
Longitudinal spacing models describe the speed--spacing states that vehicles maintain, while longitudinal controllers regulate vehicle motion through spacing and speed errors. Less explicit is how one spacing state should be connected to another through a finite, physically meaningful vehicle transition. This paper develops such a connection through $r$-Safety and Adaptive Kinematic Smoothing (AKS). The $r$-Safety relation represents speed-dependent spacing through a braking-normalized reserve, whose change, together with leader motion, determines the follower displacement required during a transition. AKS then converts this displacement into a low-dimensional analytical trajectory with consistent position, speed, spacing, and acceleration.

Experiments using controlled platoons and NGSIM trajectories identify $r$-Safety values of $0.307$ and $0.249$--$0.281$, respectively, and show that AKS can both represent observed transitions and generate prospective references under prescribed endpoint states, transition horizons, and leader motion. In same-controller ablations, adding AKS reduces integrated acceleration effort by $21$--$47\%$ and integrated squared jerk by $73$--$81\%$ across the controlled and naturalistic datasets. FHWA ACC/CACC trajectories further provide an empirical illustration of finite-horizon execution under automated following, showing how spacing adjustment can continue beyond the prescribed speed transition. The framework therefore provides a direct path from a speed--spacing relation to an analytical finite transition and its use as a longitudinal control reference.
\end{abstract}

\noindent\textbf{Keywords:} car following; speed--spacing relation; $r$-Safety;
Adaptive Kinematic Smoothing; finite transition; longitudinal control

\bigskip

\section{Introduction}
\label{sec:introduction}

Longitudinal spacing lies at the intersection of safety and efficiency in transportation systems. Additional separation provides room for perception, response, braking, and disturbance accommodation, whereas tighter spacing improves spatial utilization and traffic capacity. How spacing grows with speed is therefore fundamental to both microscopic vehicle interaction and system-level performance. For automated vehicles, this tradeoff becomes a design question about how spacing should be maintained and adjusted as traffic conditions change.

A defining strength of classical car-following theory is its ability to express vehicle interaction through a small number of interpretable quantities. The linear spacing rule of \citet{pipes1953operational} and the space--time translations of \citet{newell2002simplified} provide particularly economical descriptions of following behavior. In steady following, a linear relation combines standstill clearance with a time-gap component. Braking physics introduces a different scale because stopping distance grows with the square of speed at a fixed deceleration. This reasoning underlies braking-based spacing and safe-speed models \citep{herrey1945principles,gipps1981behavioural}. The additional spacing maintained at higher speed can thus be interpreted as a reserve that is released as speed falls and replenished as speed rises.

The transition between operating speeds brings this spacing reserve into the control problem. Adaptive cruise control (ACC) and cooperative adaptive cruise control (CACC) regulate following motion through feedback, with communication enabling coordinated responses within a platoon. Representative studies show how controller design governs disturbance propagation \citep{ploeg2014controller} and how following-gap choices affect freeway capacity \citep{shladover2012impacts}. Meanwhile, analytical trajectory methods construct vehicle motion from boundary conditions and physical constraints \citep{zhou2017parsimonious,malikopoulos2019closedform}. Together, these developments motivate a direct question: how can the spacing associated with a change in speed guide the motion that realizes that change?

A spacing relation specifies the separation associated with each operating speed. Its change therefore specifies how much distance the follower must gain or yield relative to the leader. Yet this relative-displacement requirement leaves open when the speed adjustment should occur and how the spacing reserve should evolve during the maneuver. Different trajectories can reach the same speed while producing different terminal gaps, acceleration demands, and intermediate clearances. The methodological opportunity is to carry the physical meaning of spacing through the entire transition, from the reserve released or replenished, to the motion required to reach the next state, and finally to a finite-horizon reference that can be executed by a car-following controller.

This paper addresses three research questions. Here, a speed--spacing state refers to a following speed and its associated steady or quasi-steady spacing, and a transition is the vehicle motion connecting two such states.
\begin{enumerate}
    \item \textbf{Spacing representation:} Can the relationship between following speed and spacing be represented parsimoniously in a form that is both physically interpretable and empirically meaningful?
    \item \textbf{Transition realization:} How can a change between two speed--spacing states be realized as an analytical finite transition that satisfies the required follower displacement?
    \item \textbf{Reference generation and execution:}  Can such an analytical transition serve as a motion reference for longitudinal control, and how does this guidance affect the realization of the desired spacing transition?
\end{enumerate}

The proposed framework combines $r$-Safety with Adaptive Kinematic Smoothing (AKS). $r$-Safety represents the nonlinear component of spacing on a braking-distance scale, making its release and replenishment explicit. The spacing change, together with leader motion, determines the follower displacement required between two states. AKS realizes this displacement through a low-dimensional analytical transition family and provides consistent position, speed, and acceleration references for longitudinal execution.

The paper makes three contributions.
\begin{enumerate}
    \item It develops $r$-Safety as a parsimonious spacing relation that combines the classical Newell linear reference with a braking-normalized reserve. The resulting formulation gives a physically interpretable representation of nonlinear spacing and links changes in spacing to the follower displacement required between operating states.
    \item It develops Adaptive Kinematic Smoothing as a low-dimensional analytical family for finite car-following transitions. The KS1, KS1A, and KS2 structures translate a required spatial adjustment into explicit vehicle motion while preserving endpoint consistency and admissible longitudinal behavior.
    \item It establishes an analytical bridge from spacing state to longitudinal execution. The analytical transition provides a finite-horizon motion reference that can be passed to a feedback controller, allowing the effect of explicit transition guidance to be evaluated
    separately from the underlying control law.
\end{enumerate}

The empirical analysis follows this progression from spacing representation to
transition realization and execution. Controlled platoon experiments
\citep{jiang2015experimental,zheng2023experimental} support identification of
the $r$-Safety relation and analysis of observed finite transitions. NGSIM
I-80 trajectories \citep{usdot2016ngsim} examine whether the same formulation
extends to naturalistic traffic conditions. These two datasets are then used for closed-loop replay of planned transitions, while FHWA ACC and CACC platooning experiments
\citep{tiernan2017cacc} provide an empirical automated-following setting for
illustrating finite-horizon spacing realization.

The remainder of the paper is organized as follows.
Section~\ref{sec:literature} reviews classical spacing representations,
finite-transition descriptions, and related longitudinal-control approaches.
Section~\ref{sec:method} develops the $r$-Safety spacing formulation and the
AKS transition construction.
Section~\ref{sec:empirical} evaluates spacing identification, transition
reconstruction, and prospective planning under controlled and naturalistic
traffic.
Section~\ref{sec:execution} examines closed-loop execution of the resulting
transition references and illustrates their application alongside empirical ACC/CACC responses.
Section~\ref{sec:conclusion} concludes the paper.

\section{Related Literature and Research Positioning}
\label{sec:literature}

\subsection{Reference states in classical traffic-flow and car-following theory}
\label{sec:lit_states}

Classical traffic-flow theory provides a low-dimensional representation of
traffic states through relationships among flow, density, and speed. In the
Lighthill--Whitham--Richards (LWR) theory, a prescribed flow--density relation
defines traffic states and governs the propagation of changes between them
\citep{lighthill1955kinematic,richards1956shock}. Daganzo later used this
framework to connect kinematic-wave descriptions with vehicle-level
interpretations of traffic motion \citep{daganzo1997fundamentals}. These
representations are deliberately parsimonious: much of the evolution of a
traffic stream can be described once a small set of state variables and their
relations have been specified.

At the vehicle level, car-following theory provides a corresponding description
in terms of speed and spacing. Pipes' rule relates desired following distance
linearly to speed \citep{pipes1953operational}. Newell's simplified
car-following model gives an especially economical representation by
translating the leader trajectory in space and time
\citep{newell2002simplified}. Under steady following, the translation produces
a linear spacing--speed relation characterized by a spatial offset and a time
shift. The resulting connection between microscopic following behavior and
macroscopic traffic states is one reason such low-order representations remain
central to traffic-flow theory \citep{wageningen2015genealogy}.

Other classical models introduce nonlinear spacing through different physical or behavioral mechanisms. Kometani and Sasaki derived a nonlinear spacing--speed relation from vehicle-motion and safety considerations \citep{kometani1961dynamic}. Gipps incorporated reaction and braking requirements into a safe-speed formulation \citep{gipps1981behavioural}, while Newell's earlier nonlinear model described delayed response through a nonlinear relation between speed and following distance \citep{newell1961nonlinear}. The Intelligent Driver Model combines desired spacing with relative-speed effects and free acceleration \citep{treiber2000congested}. These models differ substantially in interpretation and dynamics, but they share a useful role for the present discussion: each provides a mechanism by which an operating condition can be associated with a reference following state. Recent reviews further show that contemporary car-following research spans parsimonious physical models, data-driven formulations, and AV/CAV-oriented control models, while retaining the need for interpretable state and response representations \citep{xie2026cfreview}.

The natural starting point for a finite transition can therefore be viewed as
two reference states, an initial state $(v_1,s_1)$ and a target state
$(v_2,s_2)$. Their spacing values may arise from a theoretical equilibrium
relation, empirical quasi-steady observations, or a commanded spacing policy.
Specifying these states, however, does not yet specify the vehicle motion
between them.

\subsection{From reference states to finite car-following transitions}
\label{sec:lit_finite_transition}

The distinction between reference states and the motion connecting them is
already visible in classical traffic-flow theory. In an idealized
kinematic-wave representation, a vehicle may encounter a sharp change between
traffic states as it crosses a shock. Daganzo contrasted this construction with
a car-following description in which the vehicle passes through a finite
transition region: the LWR trajectory changes state at the shock, whereas the
car-following trajectory evolves continuously over a nonzero ``shock width''
\citep{daganzo1997fundamentals}. The comparison exposes a fundamental
difference between describing traffic states and describing the kinematics of
an individual vehicle as it moves between them.

Finite vehicle response also has a long experimental foundation. The General
Motors car-following studies measured how a follower responds to changes in the
motion of its leader and related acceleration to delayed relative-speed
information \citep{chandler1958traffic}. Subsequent work examined the
propagation and stability of such responses along a vehicle string
\citep{herman1959stability}, and nonlinear follow-the-leader formulations
extended the response law to depend on vehicle speed and spacing
\citep{gazis1961nonlinear}. These studies established that adjustment between
approximately steady operating conditions is not instantaneous: it involves a
finite response interval during which speed, relative motion, acceleration, and
vehicle position evolve together.

More recent experimental studies have continued this line of inquiry under
controlled vehicle-following conditions. Instrumented platoon experiments have
used prescribed changes in leader speed to distinguish approximately settled
following states from the acceleration and deceleration responses that connect
them \citep{jiang2015experimental,zheng2023experimental}. Such observations
provide direct evidence of finite transition duration, displacement, and
relative-motion evolution, while also showing that transient trajectories need
not collapse onto a single steady spacing--speed relation. In this sense,
modern experiments complement the classical car-following literature by making
the reference states and the intervening transition observable within the same
experimental setting~\citep{zhou2025stochastic}.

This distinction is central to the finite-transition problem. Even when the
initial speed, final speed, and transition duration are prescribed, many smooth
speed profiles can connect the same endpoint speeds. Those profiles need not
produce the same traveled distance, terminal spacing, peak acceleration, or
intermediate clearance. A smooth change in speed is therefore not, by itself, a
complete state-to-state transition. The spatial requirement implied by the
reference states must also be satisfied.

\subsection{Finite-transition construction under kinematic and safety constraints}
\label{sec:lit_transition}

Once the initial and target states have been specified, the intervening motion
can be treated as a constrained trajectory-generation problem. In its most
general form, the unknown is a continuous vehicle trajectory whose position,
speed, and acceleration must satisfy endpoint conditions together with
kinematic and pathwise constraints. The problem can therefore be approached
through optimal control, trajectory optimization, or structured analytical
constructions, depending on the application and the degrees of freedom retained
in the trajectory representation.

Several transportation studies have developed compact constructions for such
motion. \citet{delpiano2015kinematic} introduced finite deceleration into a
parsimonious car-following approximation. For connected automated traffic,
shooting methods construct piecewise analytical trajectories subject to
arrival, roadway, vehicle, and car-following constraints
\citep{zhou2017parsimonious,ma2017parsimonious}. \citet{wei2017dynamic}
combined simplified car-following dynamics with dynamic programming for
multi-vehicle longitudinal trajectory optimization and adaptive platoon
formation. Related optimal-control formulations derive closed-form trajectories
subject to state, control, and rear-end safety constraints
\citep{malikopoulos2019closedform}, while parsimonious multi-phase
constructions have been developed for anticipated leader braking
\citep{jin2025projected}. Together, these studies show that finite vehicle
motion can be represented through either low-dimensional analytical structures
or more general numerical trajectory optimization.

Safety and longitudinal-control methods address complementary aspects of the same motion. Responsibility-Sensitive Safety derives admissible longitudinal distance
conditions from assumptions on response and braking \citep{shalevshwartz2017rss}. Situation-aware collision-avoidance formulations have likewise adapted longitudinal spacing requirements to the prevailing driving state and perception uncertainty \citep{li2018situation}. Control barrier functions formulate safety in terms of forward invariance of an admissible set and can modify a nominal control action when that set is approached \citep{ames2017cbf}; this framework has also been applied to connected cruise control \citep{molnar2023safety, ZHANG2026106982}. ACC and CACC methods, meanwhile, regulate the execution of longitudinal motion through feedback and, in cooperative systems, communicated information \citep{ploeg2014controller,li2015overview, gong2019dynamic,vegamoor2022lossy}. Recent work has also embedded ACC feedback dynamics directly into macroscopic traffic-flow formulations, showing that transient wave propagation depends not only on equilibrium spacing but also on the underlying controller parameters \citep{li2026accwave}.

These bodies of work act on related but distinct mathematical objects.
Trajectory optimization determines a motion after boundary conditions and an
objective have been specified. Safe-distance and safe-set formulations define
conditions under which motion remains admissible. Longitudinal controllers
determine how a prescribed spacing target or motion reference is executed.
The finite-transition representation considered here lies between the
specification of the reference states and these downstream execution
mechanisms: its purpose is to construct a compact state-to-state motion that is
consistent with the required displacement and with the physical constraints
imposed on the transition.

\subsection{The finite-transition problem}
\label{sec:lit_positioning}

The preceding literature identifies three complementary objects: reference
following states, finite vehicle response between those states, and methods for
constructing, constraining, or executing vehicle motion. The question considered
here lies between the first two. Once the initial and target states are
specified, what finite vehicle motion connects them while remaining consistent
with the required displacement and kinematic constraints?

\begin{figure}[!htbp]
    \centering
    \includegraphics[width=\textwidth]{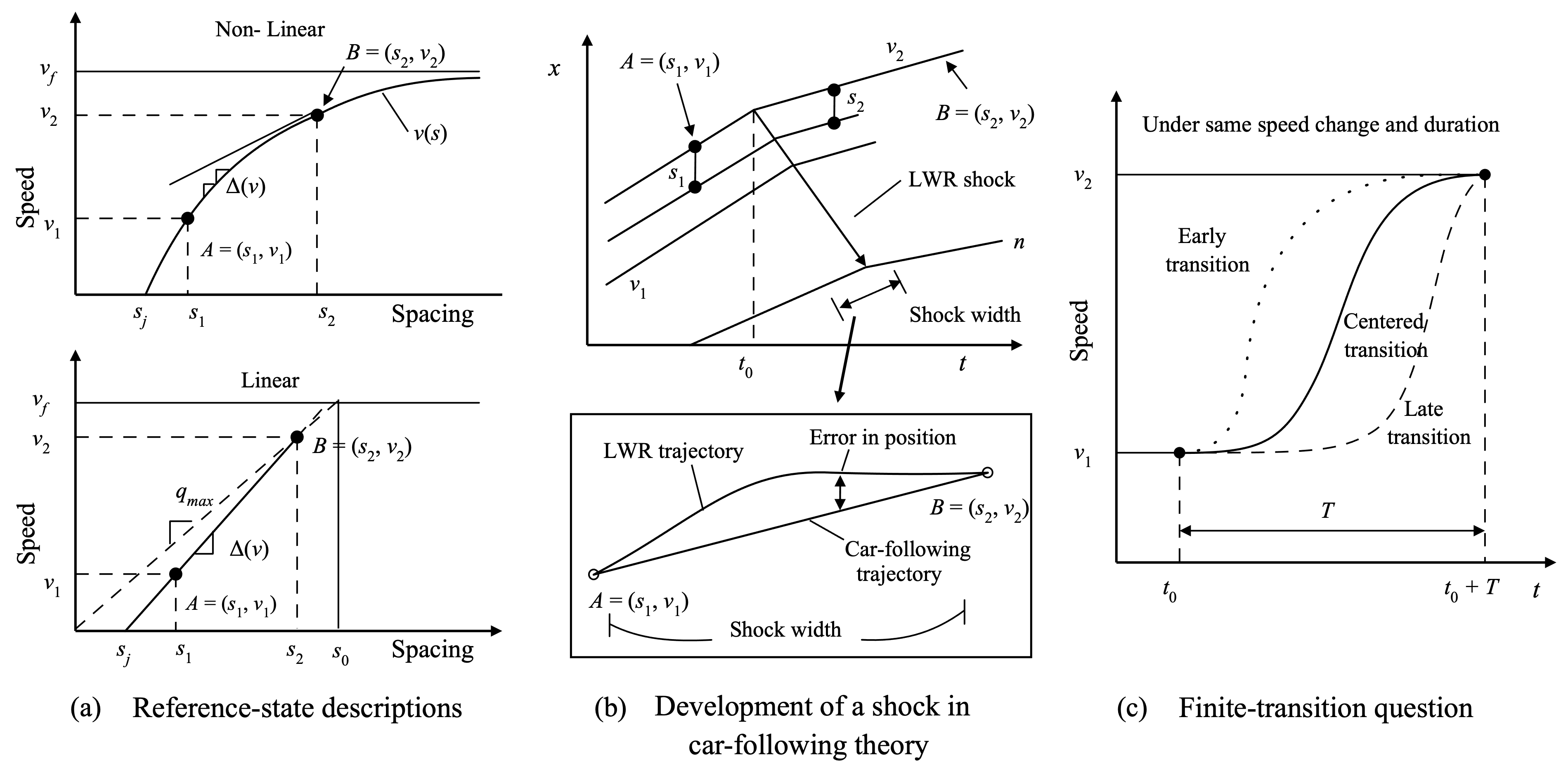}
    \caption{Classical reference-state descriptions and the finite-transition
    problem. Panels (a) and (b) are adapted from Daganzo's speed--spacing and
    shock-development constructions, with $A=(s_1,v_1)$ and $B=(s_2,v_2)$
    added as reference states. Panel (c) illustrates that identical endpoint
    speeds and transition duration can produce different follower
    displacements; reaching the prescribed target spacing therefore also
    requires displacement closure.}
    \label{fig:research_positioning}
\end{figure}

Figure~\ref{fig:research_positioning} illustrates this question using the
classical constructions of Daganzo
\citep[Fig.~4.27--4.28]{daganzo1997fundamentals}. Panel (a) shows speed--spacing
relations that define reference states, with
$A=(s_1,v_1)$ and $B=(s_2,v_2)$ added here as the initial and target states.
Panel (b) contrasts the idealized LWR state change with the finite response of
an individual vehicle. Panel (c) extends this observation to the transition
problem considered in this study.

As shown in Fig.~\ref{fig:research_positioning}(c), fixing $v_1$, $v_2$, and
$T$ does not uniquely determine the spatial outcome of the transition, because
different speed profiles generally produce different values of
$\int_0^T v(t)\,dt$. For prescribed states $A$ and $B$ and a given leader
motion, relative motion instead fixes the required follower displacement,
$\Delta x_F^{*}=\Delta x_L+s_1-s_2$. The problem studied here is therefore to
construct a parsimonious analytical transition that satisfies this displacement
closure together with the endpoint, spacing, and vehicle-motion conditions.

\section{The spacing-to-transition formulation and Adaptive Kinematic Smoothing}
\label{sec:method}

The finite-transition problem identified in Section~\ref{sec:literature}
requires two coupled elements: a physically interpretable specification of the
spacing states to be connected, and an analytical construction of the finite
motion between them. We address these two elements through $r$-Safety and
Adaptive Kinematic Smoothing (AKS), respectively.

$r$-Safety represents the speed-dependent spacing state through standstill
clearance, linear time-gap spacing, and a braking-normalized reserve. The
change between two such states, together with leader motion, determines the
road displacement that the follower must realize over the transition. AKS then
constructs analytical finite-motion candidates that realize this displacement
while connecting the prescribed endpoint speeds. The resulting motion provides
mutually consistent position, speed, and acceleration references for
longitudinal execution.

\subsection{$r$-Safety spacing and displacement formulation}
\label{sec:boundary_displacement}

A finite transition between two following states requires both the endpoint
speeds and the spacings associated with those speeds. We therefore begin with
the physical components used to describe the endpoint spacing states.

\subsubsection{Braking-normalized spacing reserve}

Newell's trajectory-translation model provides a parsimonious linear reference,

\begin{equation}
S_{\mathrm{N}}(v)=d_0+\tau v,
\label{eq:newell_spacing}
\end{equation}
where $d_0$ is the standstill clearance and $\tau$ is the time-gap parameter
\citep{newell2002simplified}. Equation~\eqref{eq:newell_spacing} accounts for
standstill clearance and the spacing induced by a fixed time translation, but
contains no explicit braking-distance reserve.

A complementary physical scale follows from classical constant-deceleration
kinematics. For a vehicle braking from speed $v$ to rest with constant
deceleration magnitude $a_b>0$,

\begin{equation}
0=v^2-2a_b d_{\mathrm{brake}},
\qquad
d_{\mathrm{brake}}=\frac{v^2}{2a_b}.
\label{eq:braking_distance}
\end{equation}

The quadratic speed dependence in Eq.~\eqref{eq:braking_distance} is the
stopping-distance structure that also appears in classical braking-based
spacing formulations such as \citet{herrey1945principles}.

The $r$-Safety formulation combines the linear spacing reference in
Eq.~\eqref{eq:newell_spacing} with the braking-distance scale in
Eq.~\eqref{eq:braking_distance},

\begin{equation}
S_r(v)
=
d_0+\tau v
+
r\,\frac{v^2}{2a_{\mathrm{ref}}},
\label{eq:r_safety}
\end{equation}
where $a_{\mathrm{ref}}>0$ is a prescribed reference braking magnitude and
$r$ is a dimensionless reserve parameter. The cases $r=0$ and $r=1$ have
direct interpretations: $r=0$ recovers the linear Newell reference, whereas
$r=1$ adds one full reference stopping-distance reserve. Intermediate values
scale this reserve continuously.

Throughout the road-vehicle analyses, we use
$a_{\mathrm{ref}}=3.4138~\mathrm{m/s^2}$, the AASHTO
stopping-sight-distance design deceleration, as a common physical
normalization \citep{aashto2018greenbook}. This reference is fixed
independently of the actual acceleration realized during a vehicle transition.
If the empirical spacing relation is written as
$S(v)=d_0+\tau v+C v^2$, then the dimensional curvature and the normalized
reserve parameter are related by $r=2a_{\mathrm{ref}}C$.

The $r$-Safety relation therefore separates the spacing associated with a
following state into a standstill component, a linear time-gap component, and
a speed-dependent reserve. The reserve is a spacing-state quantity measured
relative to a stated braking-distance scale; it is not, by itself, a formal
collision-risk measure or collision-avoidance guarantee.

\begin{figure}[!htbp]
    \centering
    \includegraphics[width=\textwidth]{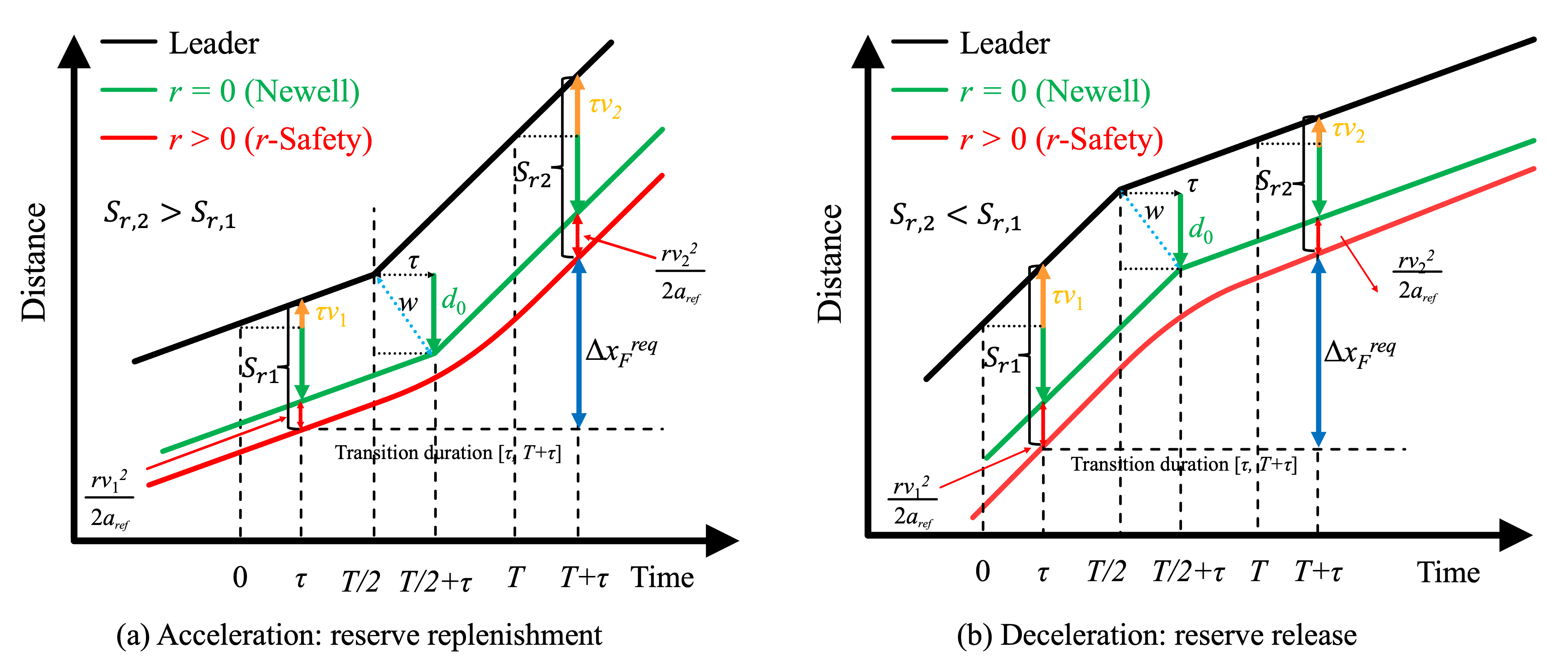}
    \caption{$r$-Safety spacing reserve relative to the linear Newell
    reference. The green trajectories denote the $r=0$ reference and the red
    trajectories include the additional reserve. (a) During acceleration, the
    reserve increases from the initial to the target state and must be
    replenished. (b) During deceleration, the reserve decreases and is
    released.}
    \label{fig:r_safety_reserve}
\end{figure}

Figure~\ref{fig:r_safety_reserve} illustrates this decomposition and its change
between two following states. The green trajectories represent the Newell
reference ($r=0$), while the red trajectories include the additional
$r$-Safety reserve. At endpoint $i$, the nonlinear reserve is
$r v_i^2/(2a_{\mathrm{ref}})$, superimposed on the linear time-gap spacing
$\tau v_i$. During acceleration, $v_2>v_1$, so the higher-speed state requires
a larger reserve: the nonlinear spacing component is replenished as the
vehicle moves from the initial to the target state. During deceleration,
$v_2<v_1$, the required reserve becomes smaller and part of the spacing
maintained at the initial state is released. The figure thus describes a
change in the spacing state itself, independently of the particular finite
trajectory used to connect the two endpoints.

\subsubsection{Conditional safety interpretation of the reserve}
\label{sec:conditional_safety}

To make the safety meaning of the braking-normalized reserve explicit, we
consider how it affects the clearance remaining after a specified
leader-braking event. Let event $e$ begin at time $t=t_e$, when the realized
net gap is $g(t_e)$ and the leader and follower speeds are $v_L(t_e)$ and
$v_F(t_e)$, respectively. Let $b_L>0$ and $b_F>0$ denote the constant braking
magnitudes of the leader and follower, and let $\delta$ denote the follower
response delay. If the leader begins braking at $t_e$ and the follower
maintains its pre-event speed during the response delay before braking, the
residual clearance after the assumed braking responses is

\begin{equation}
M_e
=
g(t_e)
+\frac{v_L^2(t_e)}{2b_L}
-v_F(t_e)\delta
-\frac{v_F^2(t_e)}{2b_F},
\label{eq:event_residual_clearance}
\end{equation}
where $M_e>0$ denotes positive residual clearance, $M_e=0$ denotes boundary
contact, and $M_e<0$ denotes trajectory overlap under the stated braking
assumptions. Equation~\eqref{eq:event_residual_clearance} follows the
delayed-response braking logic used in dynamic errorable car-following safety
analysis \citep{przybyla2015risk}.

To isolate the contribution of the $r$-Safety reserve, consider a pre-event
state on the spacing reference, such that
$g(t_e)=S_r\!\left(v_F(t_e)\right)$. Substituting
Eq.~\eqref{eq:r_safety} into Eq.~\eqref{eq:event_residual_clearance} and
holding the braking environment fixed gives

\begin{equation}
\frac{\partial M_e}{\partial r}
=
\frac{v_F^2(t_e)}{2a_{\mathrm{ref}}}
>0.
\label{eq:r_clearance_sensitivity}
\end{equation}

Thus, under the same pre-event speed and braking conditions, increasing $r$
adds a directly quantifiable amount of residual clearance. In particular, a
unit increase in $r$ contributes one reference stopping-distance unit,
$v_F^2(t_e)/(2a_{\mathrm{ref}})$, to the post-braking clearance. This provides
the conditional safety interpretation of the $r$-Safety reserve. The actual
event outcome still depends on response delay, braking capabilities, vehicle
states, and realized spacing, so $r$ itself should not be interpreted as a
crash probability.

Across drivers, controllers, and traffic conditions, these event quantities
may vary and can ultimately be treated probabilistically. Time-varying
car-following behavior and random-parameter safety models provide possible
routes for representing such heterogeneity and sparse adverse outcomes
\citep{taylor2015heterogeneity,pantangi2021crosswalks}. The present study
retains the deterministic reserve and finite-transition interpretation, while
a probabilistic safety analysis is left for future work.

\subsubsection{Spacing-to-displacement bridge}

The change in spacing between the two endpoint states has a direct kinematic
consequence. Let $x_L(t)$ and $x_F(t)$ denote the front-bumper positions of
the leader and follower, and let $L_L$ denote the leader length. The net gap
is $g(t)=x_L(t)-x_F(t)-L_L$. Consider a transition of duration $T$ from
$A=(s_1,v_1)$ to $B=(s_2,v_2)$, where $s_1$ and $s_2$ are the initial and
target net gaps. If $\Delta x_L$ and $\Delta x_F$ denote the corresponding
leader and follower displacements, relative motion requires

\begin{equation}
\Delta x_F^{\mathrm{req}}
=
\Delta x_L+s_1-s_2.
\label{eq:displacement_bridge}
\end{equation}

Equation~\eqref{eq:displacement_bridge} is the spacing-to-motion bridge of the
framework. A reduction in spacing requires the follower to gain distance
relative to the leader, whereas an increase in spacing requires it to yield
distance. Once the two spacing states and leader displacement are specified,
the total follower displacement is therefore fixed even though the temporal
profile by which that displacement is achieved remains undetermined.

For a fixed $r$-Safety relation, an observed transition need not begin exactly
on the population-level fitted curve. The construction therefore preserves the
initial offset from the fitted relation and uses $r$-Safety to determine the
change in reference spacing between the endpoint speeds. The required
follower displacement can then be written directly as

\begin{equation}
\begin{aligned}
\Delta x_F^{\mathrm{req}}
&=
\Delta x_L-
\left[S_r(v_2)-S_r(v_1)\right] \\[1mm]
&=
\Delta x_L
-\tau(v_2-v_1)
-r\,\frac{v_2^2-v_1^2}{2a_{\mathrm{ref}}}.
\end{aligned}
\label{eq:r_displacement_bridge}
\end{equation}

The two terms in the spacing change have distinct interpretations. The linear
term $\tau(v_2-v_1)$ is the adjustment associated with the time-gap component,
whereas $r(v_2^2-v_1^2)/(2a_{\mathrm{ref}})$ is the change in
braking-normalized reserve. During acceleration, both the time-gap spacing and
the nonlinear reserve increase, so the follower must travel less relative to
the leader in order to establish the larger target spacing. During
deceleration, the spacing requirement contracts and previously maintained
reserve is released, allowing additional relative advance by the follower.
Thus, $r$-Safety not only specifies the endpoint spacing states but also
determines how a speed change modifies the road displacement required of the
follower.

The same relative-motion formulation extends naturally to transitions between
different following regimes, for which the initial and target states may be
governed by different spacing relations and different $r$-Safety reserves.
Appendix~\ref{app:cross_regime} gives this cross-regime formulation and its
equal-speed special case. The equilibrium micro-to-macro implications of the
$r$-Safety spacing relation are summarized separately in
Appendix~\ref{app:micro_macro}.
\subsection{Analytical transition family}
\label{sec:aks_family}

The spacing-to-displacement formulation determines the total follower
displacement that must be realized over the transition, but does not determine
how that displacement is distributed in time. Adaptive Kinematic Smoothing
(AKS) represents this remaining finite-motion problem through a low-dimensional
family of analytical speed profiles. Let $v(t)$ denote the follower speed over
$t\in[0,T]$, with $\Delta v=v_2-v_1$. For approximately steady endpoint
states, the transition satisfies

\begin{equation}
v(0)=v_1,\qquad
v(T)=v_2,\qquad
a(0)=a(T)=0,\qquad
\int_0^T v(t)\,dt=\Delta x_F^{\mathrm{req}} .
\label{eq:aks_conditions}
\end{equation}

The first three conditions specify the endpoint motion, while the integral
condition fixes the road displacement accumulated during the transition.
The planning formulation considered here is conditioned on prescribed boundary
information: the initial and target speed--spacing states, the transition
horizon $T$, and the leader motion over $[0,T]$ are treated as inputs to the
transition construction. AKS then solves the finite-motion problem associated
with these conditions.

AKS represents this problem through three low-dimensional analytical
structures. KS1 gives a symmetric single-phase transition, KS1A allows the
timing of the same overall speed adjustment to shift within a single phase,
and KS2 introduces an intermediate state and two smooth phases.

\begin{figure}[!htbp]
    \centering
    \includegraphics[width=\textwidth]{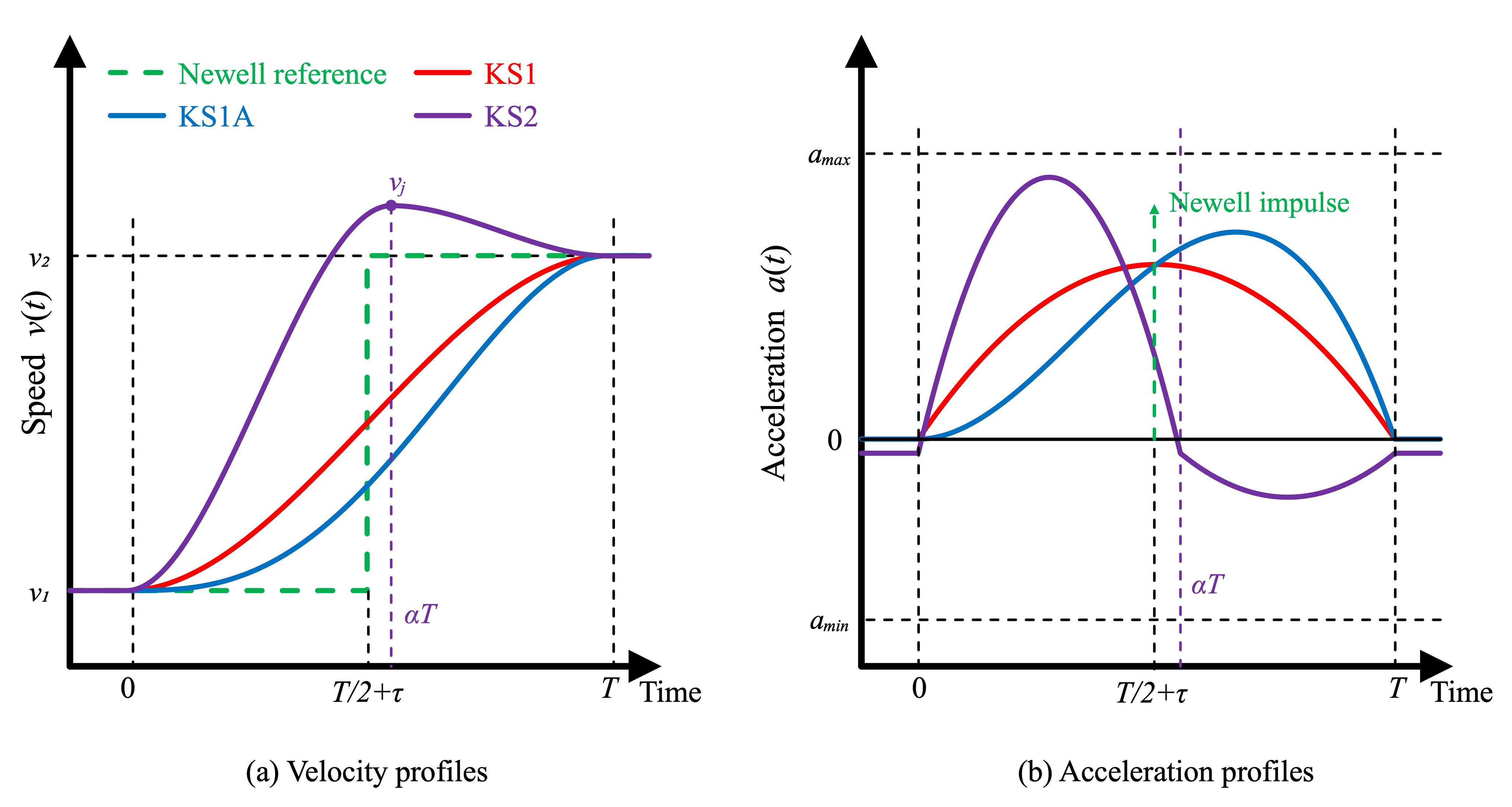}
    \caption{Representative analytical transition structures in the AKS
    family for an overall speed increase: (a) speed profiles and
    (b) acceleration profiles. KS1 is symmetric, KS1A shifts the timing of the
    speed adjustment within a single phase, and the illustrated KS2 uses two
    phases with a temporary speed overshoot.}
    \label{fig:aks_family}
\end{figure}

Figure~\ref{fig:aks_family} gives a representative comparison for an overall
speed increase. KS1 concentrates the acceleration symmetrically around the
middle of the interval. KS1A moves the acceleration activity earlier or later,
changing the area under the speed profile while preserving the same endpoint
speeds. KS2 further redistributes the transition through an intermediate
junction state; the illustrated case contains a temporary speed overshoot.

\subsubsection{KS1: symmetric single-phase transition}
\label{sec:ks1}

Let
\[
\xi=\frac{t}{T},
\qquad
h(\xi)=3\xi^2-2\xi^3 .
\]
The cubic interpolation satisfies $h(0)=0$, $h(1)=1$, and zero endpoint
slopes. The corresponding KS1 speed, acceleration, and traveled displacement
are

\begin{equation}
\begin{aligned}
v_{\mathrm{KS1}}(t)
&=
v_1+\Delta v\,h(\xi),\\
a_{\mathrm{KS1}}(t)
&=
\frac{6\Delta v}{T}\xi(1-\xi),\\
\ell_{\mathrm{KS1}}(t)
&=
v_1t+
\Delta v\,T
\left(
\xi^3-\frac{\xi^4}{2}
\right).
\end{aligned}
\label{eq:ks1_profile}
\end{equation}

The acceleration profile is symmetric about $\xi=1/2$, where its magnitude
reaches its maximum. Integrating the cubic speed profile over the complete
transition gives

\begin{equation}
\Delta x_F^{\mathrm{KS1}}
=
\frac{v_1+v_2}{2}T,
\qquad
a_{\mathrm{peak}}
=
\frac{3|\Delta v|}{2T},
\qquad
j_{\mathrm{peak}}
=
\frac{6|\Delta v|}{T^2}.
\label{eq:ks1_properties}
\end{equation}

The mean acceleration magnitude is $|\Delta v|/T$, giving the
mean-to-peak ratio $2/3$. More importantly, the KS1 displacement is completely
determined once $v_1$, $v_2$, and $T$ are specified. KS1 therefore provides a
natural symmetric reference, but it has no independent timing parameter with
which to alter the accumulated travel. If
$\Delta x_F^{\mathrm{KS1}}>\Delta x_F^{\mathrm{req}}$, the follower travels
too far relative to the prescribed terminal spacing; if it is smaller, the
follower retains additional spacing at the end of the transition.

\subsubsection{KS1A: asymmetric single-phase transition}
\label{sec:ks1a}

KS1A retains a single monotone transition while allowing the acceleration
activity to shift within the interval. This timing flexibility is introduced
through the normalized beta-function kernel

\begin{equation}
\begin{aligned}
g_{p,q}(\xi)
&=
\frac{\xi^p(1-\xi)^q}
{B(p+1,q+1)},\\
a_{\mathrm{KS1A}}(t)
&=
\frac{\Delta v}{T}\,g_{p,q}(\xi),\\
v_{\mathrm{KS1A}}(t)
&=
v_1+
\Delta v\,I_{\xi}(p+1,q+1),
\end{aligned}
\label{eq:ks1a_profile}
\end{equation}
where $B(\cdot,\cdot)$ and $I_{\xi}(\cdot,\cdot)$ denote the beta function
and regularized incomplete beta function, respectively. The normalization
ensures that the acceleration integrates to the prescribed speed change
$\Delta v$. For $p,q\geq1$, acceleration vanishes at both endpoints, while
the relative values of $p$ and $q$ determine where the acceleration activity
is concentrated within the interval.

Integrating the speed gives

\begin{equation}
\begin{aligned}
\ell_{\mathrm{KS1A}}(t)
&=
v_1t+
\Delta v\,T
\left[
\xi I_{\xi}(p+1,q+1)
-
\frac{p+1}{p+q+2}
I_{\xi}(p+2,q+1)
\right],\\
\Delta x_F^{\mathrm{KS1A}}
&=
T\left[
v_1+
\Delta v\,
\frac{q+1}{p+q+2}
\right].
\end{aligned}
\label{eq:ks1a_displacement}
\end{equation}

For $\Delta v\neq0$, the required displacement can be expressed through the
normalized mean-speed coordinate

\begin{equation}
\rho
=
\frac{\Delta x_F^{\mathrm{req}}/T-v_1}{\Delta v},
\qquad
p=\kappa(1-\rho)-1,
\qquad
q=\kappa\rho-1,
\label{eq:ks1a_parameterization}
\end{equation}

with

\[
0<\rho<1,
\qquad
\kappa\geq
\max\left\{
\frac{2}{\rho},
\frac{2}{1-\rho}
\right\}.
\]

Here $\rho$ locates the required mean transition speed between the two
endpoint speeds, while $\kappa$ controls the concentration of the acceleration
profile. For an acceleration transition, a larger $\rho$ corresponds to more
of the speed increase occurring earlier and therefore to greater accumulated
travel; a smaller $\rho$ shifts the adjustment later. The direction of this
travel effect reverses for deceleration. The acceleration magnitude reaches
its peak at $\xi=p/(p+q)$.

The symmetric KS1 profile is recovered at $\rho=1/2$ and $\kappa=4$.
The condition $0<\rho<1$ also has a direct kinematic meaning: the required
mean speed $\Delta x_F^{\mathrm{req}}/T$ lies between the two endpoint speeds.
KS1A therefore represents the range of displacement adjustments that can be
achieved by redistributing a monotone single-phase speed change in time.

\subsubsection{KS2: two-phase transition}
\label{sec:ks2}

KS2 introduces an intermediate junction state and divides the transition into
two smooth phases. Let $\alpha\in(0,1)$ denote the phase split,

\[
T_1=\alpha T,
\qquad
T_2=(1-\alpha)T,
\]
and let $v_J$ denote the junction speed. The first phase connects $v_1$ to
$v_J$, and the second connects $v_J$ to $v_2$. Using the same cubic kernel
$h(\cdot)$ within each phase gives

\begin{equation}
v_{\mathrm{KS2}}(t)
=
\begin{cases}
v_1+(v_J-v_1)
h\!\left(\dfrac{t}{T_1}\right),
&
0\leq t\leq T_1,
\\[3mm]
v_J+(v_2-v_J)
h\!\left(\dfrac{t-T_1}{T_2}\right),
&
T_1<t\leq T .
\end{cases}
\label{eq:ks2_profile}
\end{equation}

Each phase has zero acceleration at its endpoints. The two segments therefore
meet continuously in both speed and acceleration at $t=\alpha T$.
Integrating the two phases gives

\begin{equation}
\Delta x_F^{\mathrm{KS2}}
=
\frac{T_1}{2}(v_1+v_J)
+
\frac{T_2}{2}(v_J+v_2),
\qquad
v_J(\alpha)
=
\frac{2\Delta x_F^{\mathrm{req}}}{T}
-
\alpha v_1
-
(1-\alpha)v_2 .
\label{eq:ks2_displacement}
\end{equation}

Equation~\eqref{eq:ks2_displacement} gives a useful interpretation of the two
KS2 parameters. The phase split $\alpha$ determines when the transition
changes from the first cubic segment to the second, while displacement closure
determines the corresponding junction speed $v_J$. Varying $\alpha$ therefore
changes both the timing of the transition and the intermediate state through
which the required displacement is realized.

\begin{figure}[!htbp]
    \centering
    \includegraphics[width=\textwidth]{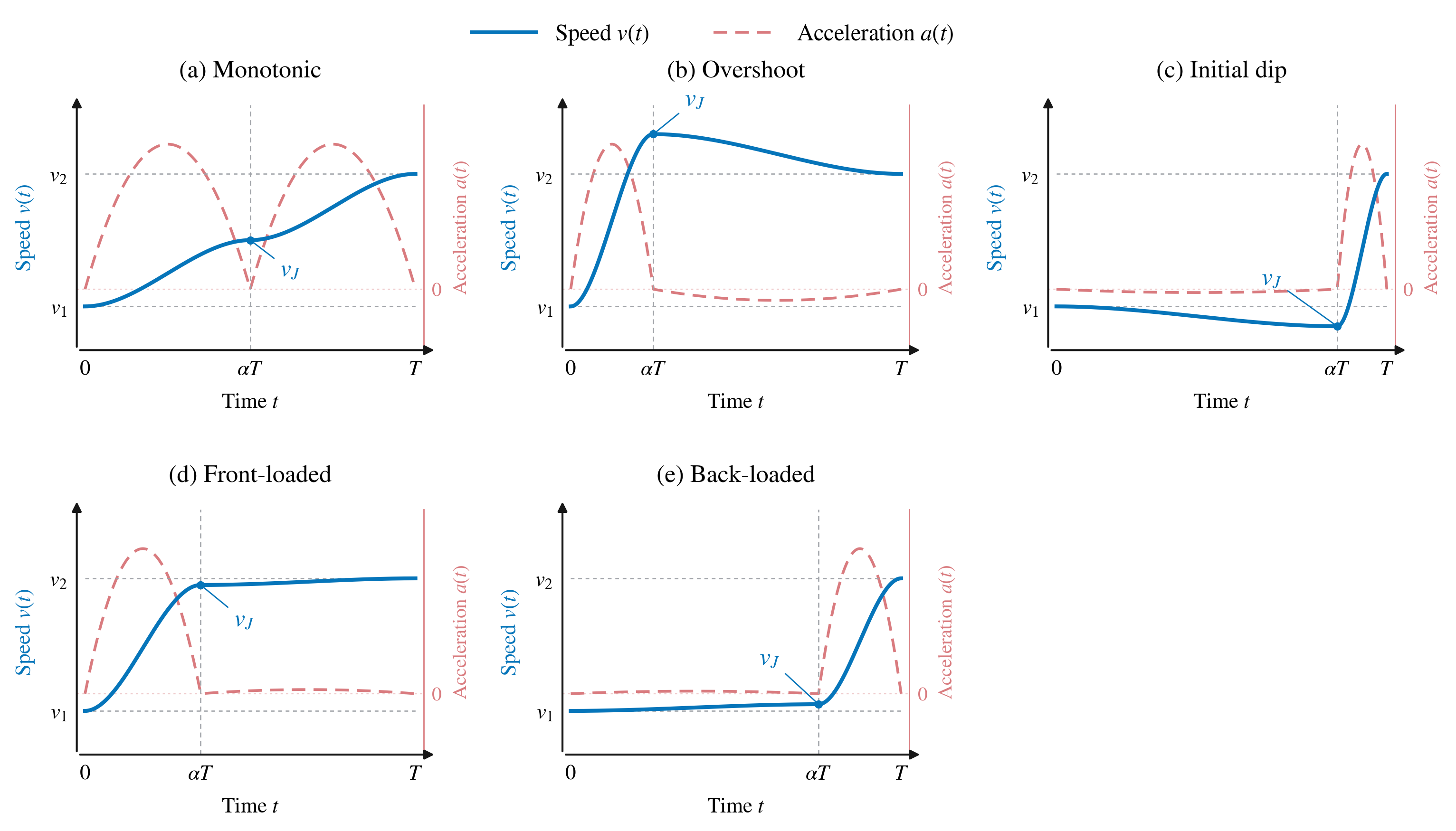}
    \caption{Representative KS2 structures generated by the junction speed
    and phase timing: (a) monotone, (b) overshoot, (c) initial dip,
    (d) front-loaded, and (e) back-loaded. Solid curves denote speed and
    dashed curves denote acceleration.}
    \label{fig:ks2_shapes}
\end{figure}

The location of $v_J$ relative to the two endpoint speeds gives several
distinct transition shapes. Figure~\ref{fig:ks2_shapes} illustrates these
five characteristic cases. For an overall speed increase,
$v_1<v_J<v_2$ produces a monotone two-phase transition. A junction close to
$v_2$ completes most of the speed adjustment early and gives a front-loaded
profile, whereas a junction close to $v_1$ produces a back-loaded profile.
If $v_J>v_2$, the vehicle temporarily overshoots the target speed before
returning to $v_2$; if $v_J<v_1$, it first dips below the initial speed before
accelerating toward the target. 

\subsection{Candidate admissibility and transition selection}
\label{sec:aks_selection}

The analytical structures in Section~\ref{sec:aks_family} generate candidate
motions under the prescribed endpoint conditions and displacement requirement.
The resulting candidates are screened for physical admissibility and then
selected according to a fixed motion-quality criterion.

Let
$\mathcal{K}=\{\mathrm{KS1},\mathrm{KS1A},\mathrm{KS2}\}$
denote the AKS family, and let $\theta_k$ collect the free parameters of
member $k\in\mathcal{K}$. For a candidate traveled displacement
$\ell_k(t;\theta_k)$, the associated speed, acceleration, and jerk are
$v_k(t;\theta_k)$, $a_k(t;\theta_k)$, and $j_k(t;\theta_k)$. Given a leader
trajectory $x_L(t)$, its planned net gap is

\begin{equation}
g_k(t;\theta_k)
=
x_L(t)
-
\left[
x_F(0)+\ell_k(t;\theta_k)
\right]
-
L_L .
\label{eq:candidate_gap}
\end{equation}

A candidate is admissible when it closes the prescribed displacement within
the terminal tolerance and maintains nonnegative speed and physical spacing
throughout the transition. The admissible parameter set for member $k$ is

\begin{equation}
\begin{aligned}
\Theta_k^{\mathrm{adm}}
=
\Bigl\{
\theta_k:\;&
\left|
\ell_k(T;\theta_k)-\Delta x_F^{\mathrm{req}}
\right|
\leq \varepsilon_x,\;
v_k(t;\theta_k)\geq0,\\
&
g_k(t;\theta_k)\geq0,
\qquad
0\leq t\leq T
\Bigr\}.
\end{aligned}
\label{eq:aks_admissible_set}
\end{equation}

For KS1, the endpoint speeds and duration uniquely determine the profile and
its displacement. If the KS1 candidate satisfies the admissibility conditions
in Eq.~\eqref{eq:aks_admissible_set}, it is retained directly. Otherwise,
displacement-closing KS1A and KS2 candidates are constructed in parallel and
evaluated over their respective admissible parameter sets.

The fitted $r$-Safety relation determines the speed-dependent spacing change
used to obtain $\Delta x_F^{\mathrm{req}}$. Its evolution along a candidate
trajectory can also be examined by retaining the initial spacing offset and
forming the translated reference

\begin{equation}
\widetilde{S}_r\!\left(v_k(t)\right)
=
S_r\!\left(v_k(t)\right)
+
\left[g(0)-S_r(v_1)\right],
\qquad
m_{r,k}(t)
=
g_k(t)-\widetilde{S}_r\!\left(v_k(t)\right).
\label{eq:translated_r_margin}
\end{equation}

The margin $m_{r,k}(t)$ describes the planned spacing relative to the same
$r$-Safety relation used to define the endpoint spacing change. It is retained
as a pathwise spacing diagnostic, while physical admissibility is determined
by Eq.~\eqref{eq:aks_admissible_set}.

Among the admissible displacement-closing candidates, transition selection
follows a fixed lexicographic criterion based first on integrated acceleration
effort and then on integrated squared jerk,

\begin{equation}
\begin{aligned}
J_a(k,\theta_k)
&=
\int_0^T a_k^2(t;\theta_k)\,dt,\\
J_j(k,\theta_k)
&=
\int_0^T j_k^2(t;\theta_k)\,dt,
\end{aligned}
\qquad
(k^*,\theta^*)
=
\operatorname*{lex\,arg\,min}_{
k\in\{\mathrm{KS1A},\mathrm{KS2}\},\;
\theta_k\in\Theta_k^{\mathrm{adm}}
}
\left(J_a,J_j\right).
\label{eq:aks_selection}
\end{equation}

Thus, acceleration effort provides the primary ranking, with jerk effort used to distinguish candidates with equivalent acceleration cost. Acceleration, braking, and jerk extrema are also recorded to quantify the motion demand of the selected transition and are evaluated against the prescribed vehicle limits in the subsequent experiments. Real-world acceleration and deceleration capabilities are themselves speed-dependent, with empirical vehicle-activity envelopes becoming wider or narrower over different speed ranges \citep{liu2015speedaccel,zhou2015emissions}. The fixed limits used here therefore provide a common screening reference, while speed-dependent envelopes can be incorporated when vehicle-specific capability data are available.

The selected member defines the finite-horizon reference
$\ell^*(t)$, $v^*(t)$, and $a^*(t)$. Together with the leader trajectory,
these quantities provide a consistent planned longitudinal motion that is
passed to the reference-generation and execution interface developed in the
next subsection.

\subsection{Reference generation and longitudinal execution interface}
\label{sec:reference_interface}

The selected AKS transition provides a finite-horizon kinematic description of
the follower motion. Let the selected candidate from
Section~\ref{sec:aks_selection} be denoted by
$\ell^*(t)$, $v^*(t)$, and $a^*(t)$ over $t\in[0,T]$. The corresponding
planned follower position is

\begin{equation}
x_F^*(t)=x_F(0)+\ell^*(t),
\label{eq:planned_position}
\end{equation}
and, for a given leader trajectory $x_L(t)$, the associated planned net gap is

\begin{equation}
s^*(t)=x_L(t)-x_F^*(t)-L_L.
\label{eq:planned_gap}
\end{equation}

Thus, the spacing-to-transition formulation yields a mutually consistent
reference set
$\{x_F^*(t),\,s^*(t),\,v^*(t),\,a^*(t)\}$.
The reference is finite-horizon and internally consistent by construction:
the planned displacement satisfies the prescribed spacing transition, the speed
profile connects the prescribed endpoint states, and the acceleration profile
is the time derivative of the same planned motion.

This reference-generation role is distinct from the longitudinal feedback law
used to execute the motion. A feedback controller may regulate an
instantaneous spacing-policy reference or track an explicitly planned
finite-horizon trajectory. In the present framework, AKS provides the latter:
the spacing change is first translated into a state-to-state motion reference,
which is then supplied to the longitudinal controller for execution.

For execution, define the tracking errors relative to the planned motion as

\begin{equation}
e_s(t)=s(t)-s^*(t),
\qquad
e_v(t)=v^*(t)-v_F(t),
\label{eq:tracking_errors}
\end{equation}
where $s(t)$ and $v_F(t)$ are the realized gap and follower speed. A generic
reference-tracking command can be written as

\begin{equation}
a_{\mathrm{cmd}}(t)
=
a^*(t)
+
k_s e_s(t)
+
k_v e_v(t),
\label{eq:acc_tracking_law}
\end{equation}
subject to the vehicle's actuation limits. Equation~\eqref{eq:acc_tracking_law}
represents an ACC-type interface in which the planned AKS acceleration serves
as a feedforward term and the spacing and speed errors provide feedback
correction. When communicated leader information is available, a CACC-type
interface can be formed by augmenting the same structure with leader-state
feedforward. The essential distinction is that AKS supplies an explicit finite-horizon
transition reference in addition to the feedback regulation of spacing and
speed errors.

The distinction is particularly clear when the target spacing changes while
the endpoint speeds remain equal. Suppose the leader travels at constant speed
$v_0$, and the follower is required to begin and end at the same speed,
\(v_1=v_2=v_0.\)
If the net gap changes from $s_1$ to $s_2$ over a duration $T$, then

\begin{equation}
\Delta x_F^{\mathrm{req}}
=
\Delta x_L+s_1-s_2
=
v_0T-(s_2-s_1).
\label{eq:equal_speed_gap_change}
\end{equation}

Relative to constant-speed travel over the same duration, the required excess
or deficit displacement is therefore

\begin{equation}
\int_0^T \bigl(v^*(t)-v_0\bigr)\,dt
=
-(s_2-s_1).
\label{eq:equal_speed_area_condition}
\end{equation}

Equation~\eqref{eq:equal_speed_area_condition} gives the physical meaning of
an equal-speed gap adjustment. If $s_2>s_1$, the follower must temporarily
travel more slowly than the constant-speed baseline and then recover to
$v_0$; if $s_2<s_1$, the corresponding speed excursion is upward. In both
cases, the endpoint speed is preserved.

\begin{figure}[!htbp]
    \centering
    \includegraphics[width=0.75\linewidth]{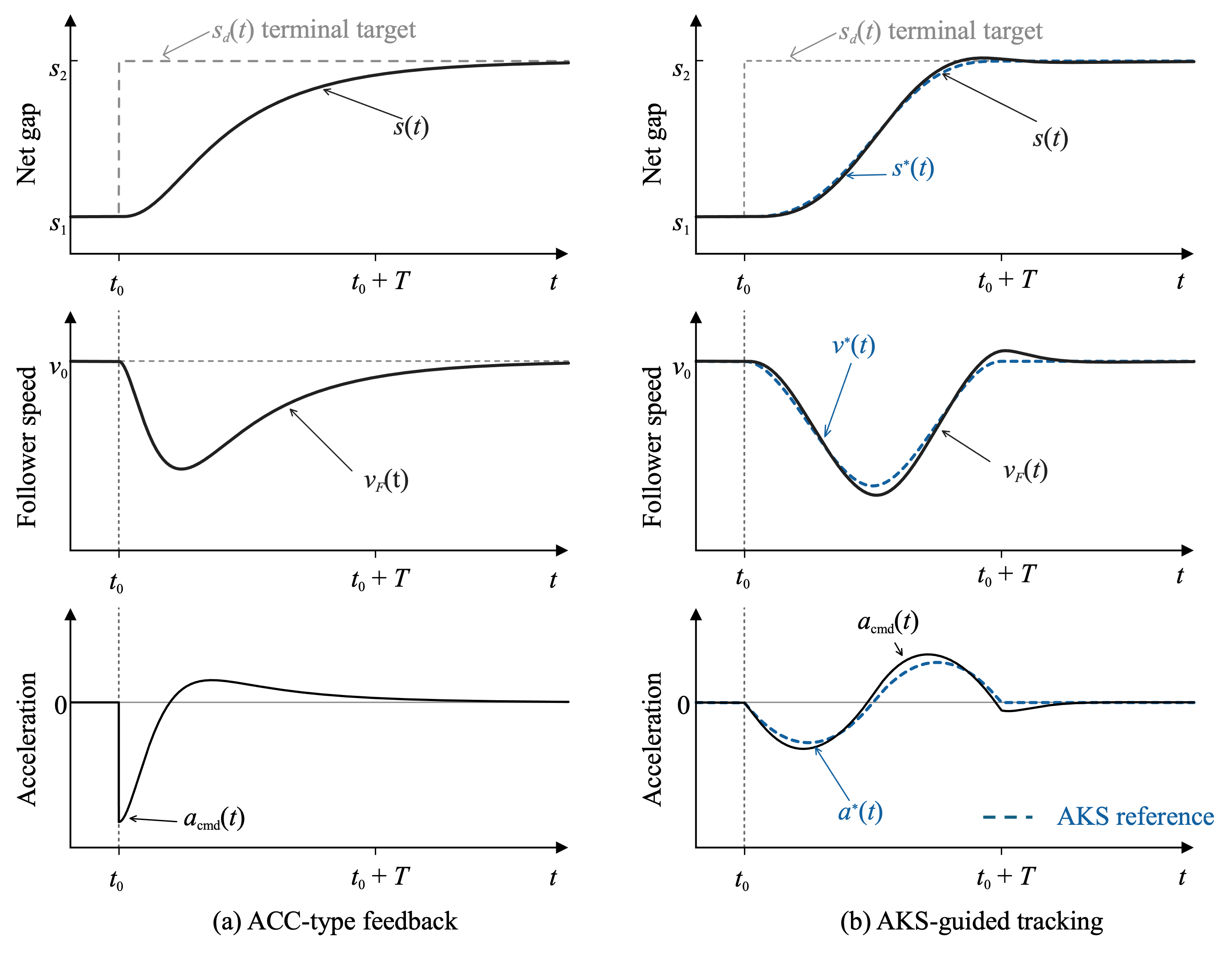}
    \caption{Feedback and AKS-guided realization of an equal-speed spacing increase.
    (a) An ACC-type feedback controller responds directly to the commanded spacing
    change from $s_1$ to $s_2$. (b) AKS first constructs a KS2 finite-horizon
    transition and supplies the resulting $s^*(t)$, $v^*(t)$, and $a^*(t)$
    references to the longitudinal controller.}
    \label{fig:reference_interface}
\end{figure}

Figure~\ref{fig:reference_interface} illustrates the case $s_2>s_1$.
The leader maintains constant speed $v_0$, while the commanded spacing changes
from $s_1$ to $s_2$ at $t=t_0$. An ACC-type feedback response acts directly
on the changed spacing target, whereas the AKS-guided case first constructs a
KS2 transition over the prescribed horizon $T$ and then tracks the resulting
spacing, speed, and acceleration references.

For online use, the endpoint states and transition horizon can be supplied by
a higher-level planner or supervisory controller, while the required
leader-motion information may be obtained from preview, prediction, or
communication. AKS then generates the corresponding finite-horizon motion
reference.

\FloatBarrier

\section{Empirical Identification and Transition Evaluation}
\label{sec:empirical}

This section evaluates the spacing-to-transition formulation using controlled
and naturalistic car-following data. The analysis proceeds from identification
of the $r$-Safety spacing relation, to reconstruction and prospective planning of
complete transitions, and finally to replication across different traffic
conditions.

\subsection{Identification of the $r$-Safety spacing reserve}
\label{sec:exp_boundary}

The controlled platoon experiments of \citet{jiang2015experimental} and
\citet{zheng2023experimental} drive a lead vehicle through a prescribed
sequence of speed plateaus and record the platoon response at 10~Hz. Five
complete runs provide 90 adjacent leader--follower pairs. The prescribed
leader-speed sequence provides both quasi-steady following states for spacing
identification and complete plateau-to-plateau transitions for AKS evaluation.

\begin{figure}[!htbp]
    \centering
    \includegraphics[width=\textwidth]{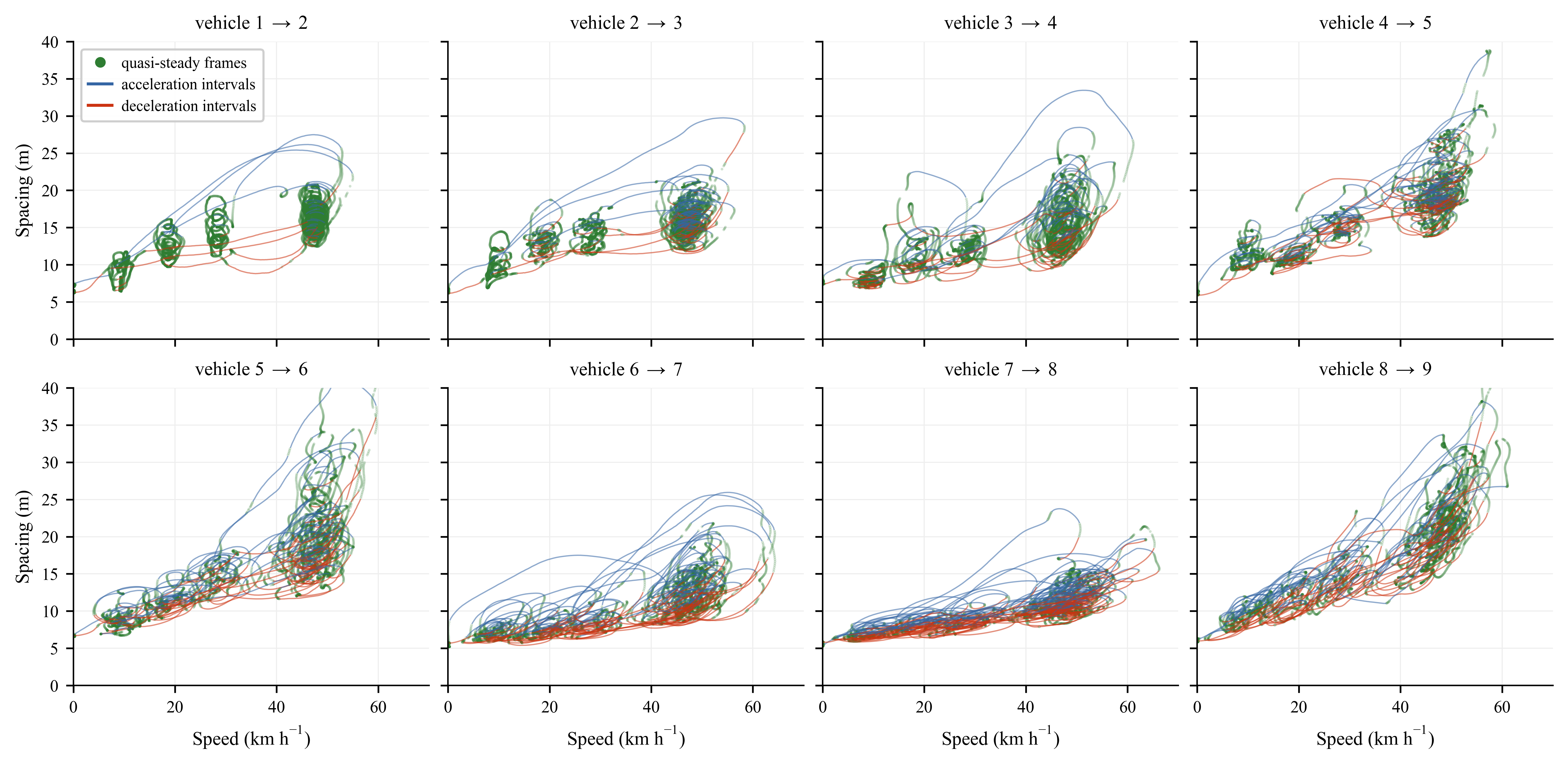}
    \caption{Speed--spacing phase diagram for the eight leader--follower pairs
    of one controlled run. Settled following forms compact clusters; the
    intervals between plateaus leave the cluster locus and follow different
    paths in acceleration and deceleration.}
    \label{fig:phase_diagram}
\end{figure}

In speed--spacing coordinates, the quasi-steady states and finite transitions
separate clearly. Figure~\ref{fig:phase_diagram} shows the eight pairs of one run. Settled
following collects into compact clusters at the commanded plateau speeds. The
intervals between clusters leave that locus entirely, tracing extended loops
that depart from the cluster pattern and return to it, with acceleration and
deceleration following visibly different paths. Steady states can therefore be
isolated by an acceleration criterion $|a_f|\leq\epsilon$, and what that
criterion excludes is structure rather than noise: the loops are the finite
transitions analyzed in Section~\ref{sec:exp_jiang}.

\begin{figure}[!htbp]
    \centering
    \includegraphics[width=\textwidth]{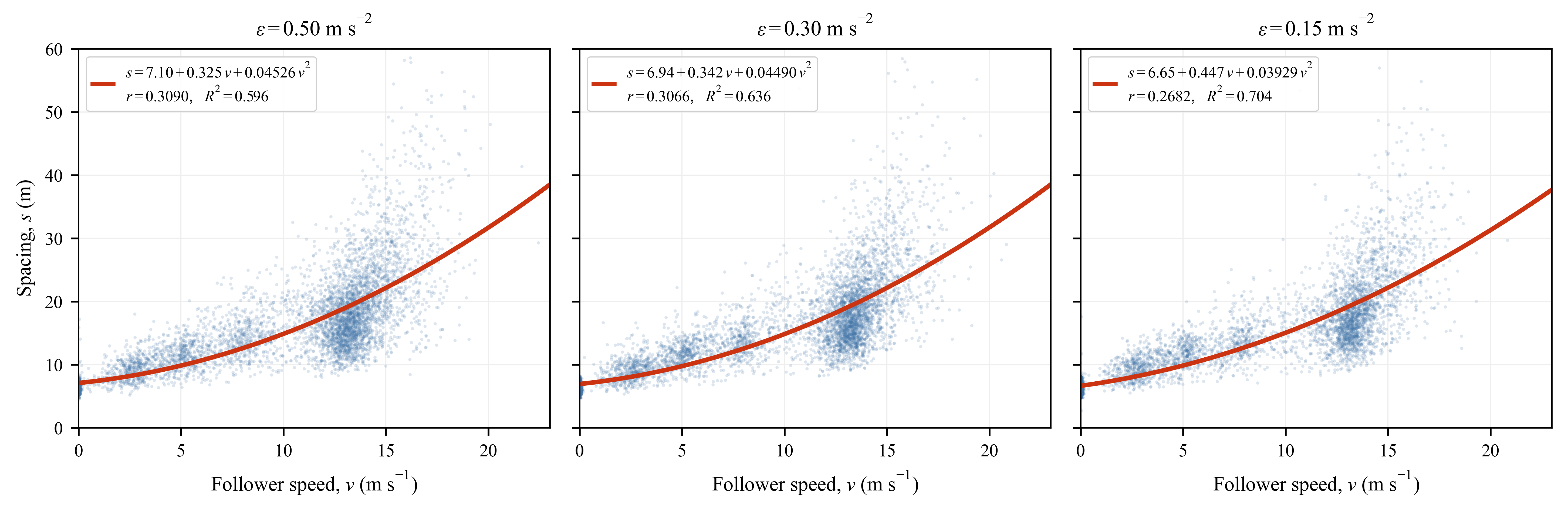}
    \caption{Quasi-steady speed--spacing relation estimated at three
    steady-state admission thresholds $\epsilon$, each with the fitted
    relation, the braking-normalized reserve $r=2a_{\mathrm{ref}}C$ and the
    coefficient of determination. Tightening $\epsilon$ reduces residual
    scatter while leaving the reserve essentially unchanged.}
    \label{fig:r_multilevel}
\end{figure}

The $r$-Safety relation in Eq.~\eqref{eq:r_safety} is estimated from the
quasi-steady frames. Figure~\ref{fig:r_multilevel} repeats the estimation at
three steady-state admission thresholds. Across the three thresholds, the
estimated reserve remains between $0.27$ and $0.31$, while the fit becomes
tighter as the steady-state criterion is made more restrictive. The
intermediate threshold $\epsilon=0.30~\mathrm{m/s^2}$ is used in the
subsequent analysis.

How large that reserve is can be assessed against the two classical limits of
the same spacing relation. At $r=0$, Eq.~\eqref{eq:r_safety} reduces to the
Newell--Pipes linear reference, whereas $r=1$ adds one complete reference
stopping-distance reserve. We compare these two limits with the fitted
$r$-Safety relation using leave-one-run-out validation. In each fold, the
spacing parameters are estimated from four controlled runs and evaluated on
the remaining run; the reported RMSE and $R^2$ are therefore calculated only
from data not used for parameter estimation.

Table~\ref{tab:reserve_levels} summarizes the out-of-sample comparison. The
fitted $r$-Safety relation gives the lowest validation error in all five
folds, with $\hat r=0.307$ and fold estimates ranging from $0.300$ to $0.320$.
The linear $r=0$ reference performs slightly worse, whereas the full
braking-distance case $r=1$ substantially overpredicts the observed spacing.
Allowing the latter to use unconstrained spacing parameters reduces its error,
but requires a nonphysical negative time gap and remains less accurate than
the linear reference. At $15~\mathrm{m/s}$, the fitted nonlinear reserve
contributes $10.1~\mathrm{m}$ of spacing, compared with $5.1~\mathrm{m}$ from
the time-gap term and $6.9~\mathrm{m}$ from standstill clearance.

\begin{table}[!htbp]
\centering
\caption{Out-of-sample validation of the fitted $r$-Safety relation against
the linear and full braking-reserve limits.}
\label{tab:reserve_levels}

\begin{tabular}{lcc}
\toprule
Spacing relation & Validation RMSE (m) & Validation $R^2$ \\
\midrule
Newell reference ($r=0$)
& 4.700 & 0.616 \\

Fitted $r$-Safety ($\hat r=0.307$)
& \textbf{4.594} & \textbf{0.633} \\

Full braking reserve ($r=1$)
& 8.654 & $-0.302$ \\
\bottomrule
\end{tabular}
\end{table}

\FloatBarrier
\subsection{Representing and planning controlled transitions with AKS}
\label{sec:exp_jiang}

The controlled experiments yield 477 complete
plateau--transition--plateau episodes, including 243 acceleration and
234 deceleration transitions, with a median duration of $26.6~\mathrm{s}$.
The same AKS family is evaluated in two modes. Reconstruction uses the
realized follower displacement and observed speed trajectory to identify the
analytical transition that best represents the maneuver. Planning withholds
the interior follower trajectory, obtains the required displacement from the
$r$-Safety spacing transition, and selects an admissible reference by $J_a$ and then
$J_j$.\footnote{Speeds are smoothed with a $2.1~\mathrm{s}$
Savitzky--Golay filter. A plateau lasts at least $5~\mathrm{s}$, with speed
standard deviation below $0.40~\mathrm{m/s}$ and absolute slope below
$0.05~\mathrm{m/s^2}$. Complete transitions satisfy
$|\Delta v|\geq1.5~\mathrm{m/s}$ and last between $1$ and $60~\mathrm{s}$.}

The controlled transitions show a clear progression in the displacement
coverage of the three analytical structures. KS1 closes 45 of the 477
realized displacements. Allowing the transition timing to shift through KS1A
increases this number to 397. For the remaining 80 transitions, the required
mean speed lies outside the interval between the endpoint speeds, so a
monotone single-phase profile cannot satisfy displacement closure. KS2 closes
all 477 transitions. Under retrospective reconstruction, KS2 provides the
best velocity representation in 419 cases, compared with 56 for KS1A and
2 for KS1.

\begin{figure}[!htbp]
    \centering
    \includegraphics[width=\textwidth]{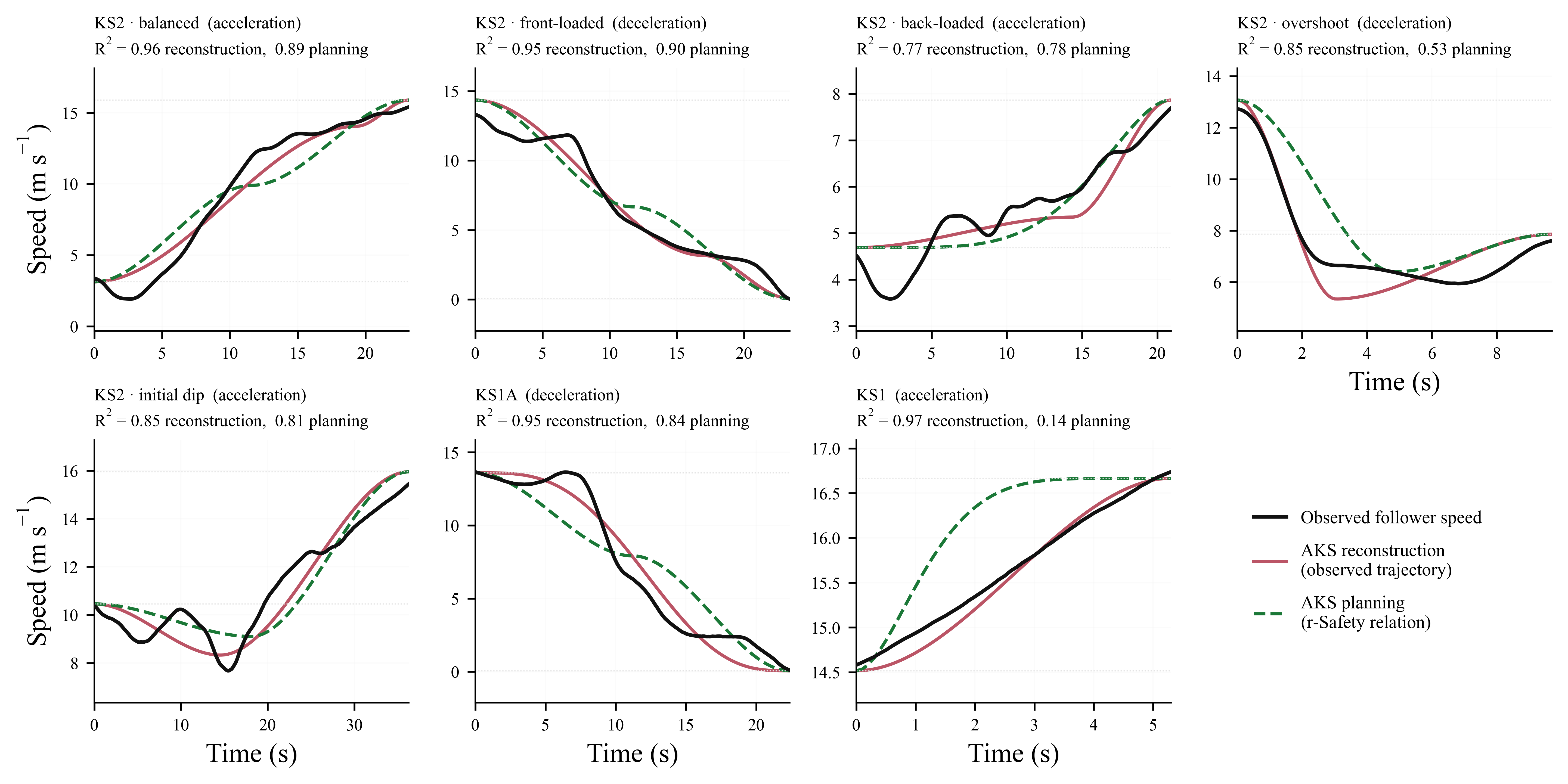}
    \caption{Representative controlled transitions across the AKS family.
    Black curves denote the observed follower speed, red curves the
    retrospective AKS reconstruction, and green dashed curves the prospective
    reference generated from the fitted $r$-Safety relation. The panels show the five
    characteristic KS2 structures together with representative KS1A and KS1
    cases. Within each structure, the displayed transition has reconstruction
    $R^2$ closest to that structure's median among events lasting at least
    $8~\mathrm{s}$ with a feasible planned reference.}
    \label{fig:transition_atlas}
\end{figure}

Figure~\ref{fig:transition_atlas} shows representative trajectory-level
examples across the analytical family. Each panel contains an observed
transition together with its retrospective AKS reconstruction and the
prospective reference generated from the fitted $r$-Safety relation. Within each
structure, the displayed case is chosen to have reconstruction performance
closest to that structure's median rather than to maximize goodness of fit.
The examples therefore illustrate typical behavior across balanced,
front-loaded, back-loaded, overshoot, initial-dip, KS1A, and KS1 transitions.

Across all 477 transitions, the selected reconstruction reaches a median
velocity $R^2$ of $0.894$. Its integrated acceleration effort also follows the
observed maneuver demand, with a transition-level Spearman correlation of
$0.72$ and a median fitted-to-observed ratio of $0.76$. The analytical
representation is substantially smoother at the jerk level, retaining only
about $15\%$ of the observed squared-jerk integral. Of the 477 best-fit
reconstructions, 437 also satisfy the physical admissibility conditions; the
40 remaining cases fail only through an intermediate negative-gap violation.

The prospective experiment then removes the observed follower displacement and
interior trajectory while retaining the prescribed boundary conditions: the
initial and target states, transition duration, and leader motion over the
horizon. Given these inputs and the fitted $r$-Safety relation, AKS constructs
the required finite-horizon reference. When compared afterward with the observed
trajectories, these references retain a median velocity $R^2$ of $0.702$.
The prescribed follower displacement differs from the realized
displacement by a median of $3.3~\mathrm{m}$, accounting for part of the
difference between reconstruction and planning performance.

Table~\ref{tab:exp_modes} summarizes the two uses of the same analytical
family. Reconstruction strongly favors the more flexible KS2 structure,
whereas planning selects KS1A and KS1 much more often. This shift is
consistent with the different objectives: reconstruction seeks the closest
representation of observed motion, while planning seeks an admissible
low-effort realization of the $r$-Safety-implied displacement.

\begin{table}[!htbp]
\centering
\caption{Summary of retrospective reconstruction and prospective planning for
the controlled transitions.}
\label{tab:exp_modes}
\small
\begin{tabular}{lcccc}
\toprule
Mode &
Displacement source &
Output &
KS1 / KS1A / KS2 &
Median $R^2$ \\
\midrule
Reconstruction &
Observed &
477 reconstructed &
2 / 56 / 419 &
0.894 \\
Planning &
Fitted $r$-Safety relation &
441 feasible references &
29 / 192 / 220 &
0.702 \\
\bottomrule
\end{tabular}
\end{table}

These results establish both roles of AKS in the controlled data: the family
provides a compact analytical representation of observed finite transitions,
and the same construction can generate prospective finite-horizon references
when the transition itself has not yet been observed.

\subsection{Generalization to naturalistic traffic}
\label{sec:exp_ngsim}

We next examine whether the same spacing-to-transition formulation extends
from controlled platoon experiments to naturalistic freeway traffic. The
NGSIM I-80 data contain three fifteen-minute periods representing congestion
buildup, the transition toward congestion, and full congestion
\citep{usdot2016ngsim}. Each period is analyzed separately using the same
$r$-Safety identification and AKS reconstruction/planning procedures as in
Sections~\ref{sec:exp_boundary} and~\ref{sec:exp_jiang}.\footnote{Passenger-car
pairs are used in all lanes. Quasi-steady frames use
$\epsilon=0.30~\mathrm{m/s^2}$, and validation folds are formed by
leader--follower pair. The three periods contain 138, 248, and 332 complete
transitions, respectively.}

\begin{figure}[!htbp]
    \centering
    \includegraphics[width=\textwidth]{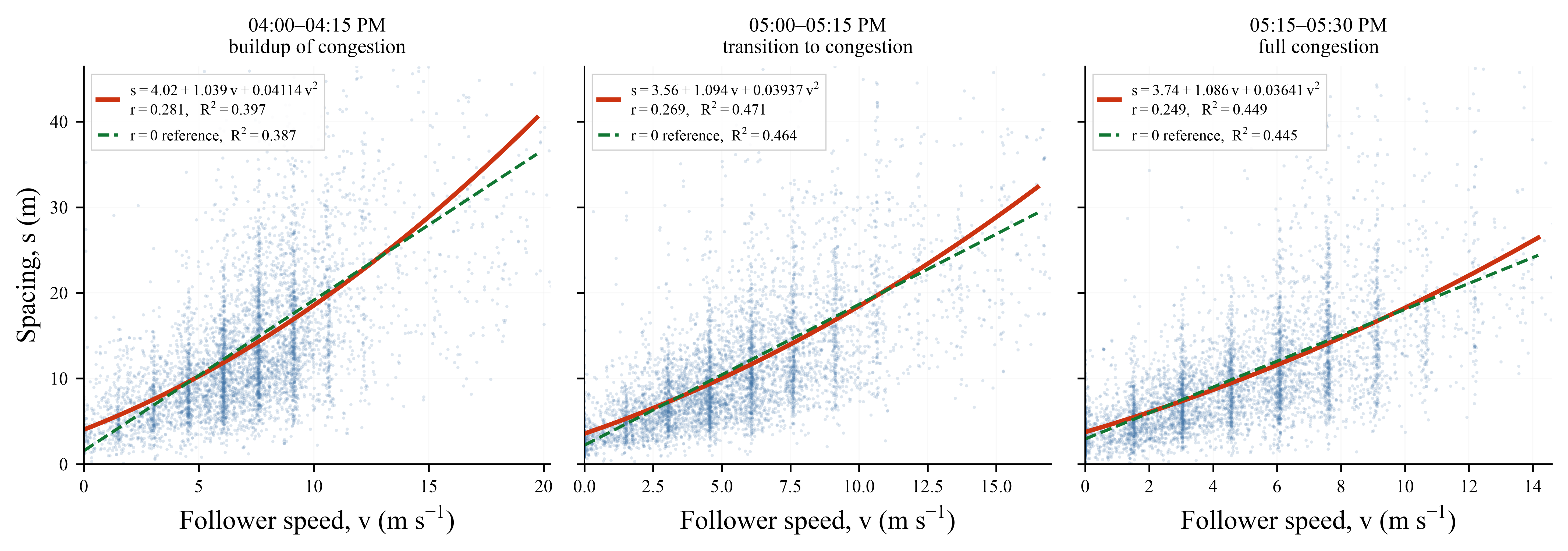}
    \caption{$r$-Safety identification for the three NGSIM I-80 periods.
    The fitted $r$-Safety relation and the linear $r=0$ reference are shown
    over the observed speed range of each period; $R^2$ is evaluated out of
    sample.}
    \label{fig:ngsim_periods}
\end{figure}

The $r$-Safety relation remains identifiable across all three traffic
conditions. Figure~\ref{fig:ngsim_periods} gives estimated reserves of
$0.281$, $0.269$, and $0.249$ from congestion buildup to full congestion,
compared with $0.307$ in the controlled platoon data. The fitted relation
outperforms the linear $r=0$ reference in each period, while the full-braking
case $r=1$ remains less accurate. The difference from the linear reference is
smaller than on the test track, particularly under congestion, where lower
speeds reduce the contribution of the nonlinear reserve.

The AKS results show a similar transfer of the transition construction.
Retrospective reconstruction achieves median velocity $R^2$ values of
$0.824$, $0.835$, and $0.774$ across the three periods. When the interior
follower trajectory is withheld, prospective planning still produces a
feasible reference for $92\%$, $88\%$, and $91\%$ of the extracted
transitions. The corresponding planning $R^2$ values are $0.477$, $0.413$,
and $0.397$, indicating that a substantial part of the realized transition structure can
still be reproduced from the prescribed endpoint states, transition horizon,
leader motion, and spacing relation.

Table~\ref{tab:transfer} summarizes the controlled and naturalistic results.
The estimated reserve decreases modestly as congestion develops, and the
agreement between planned and observed trajectories also weakens. In contrast,
the fraction of transitions admitting a feasible planned reference remains
close to $90\%$ throughout. Reconstruction likewise remains near $0.82$ in the
first two NGSIM periods before decreasing under full congestion. These results
suggest that denser traffic primarily reduces how closely a low-dimensional
analytical transition resembles the realized driving trajectory, while the
ability to construct an admissible transition is largely retained.

\begin{table}[!htbp]
\centering
\caption{Generalization of the $r$-Safety and AKS formulations from the
controlled platoon to naturalistic NGSIM traffic.}
\label{tab:transfer}

\small
\setlength{\tabcolsep}{7pt}
\renewcommand{\arraystretch}{1.08}

\begin{tabular}{l c @{\hskip 12pt} ccc}
\toprule
&
\multicolumn{1}{c}{Controlled}
&
\multicolumn{3}{c}{NGSIM I-80}
\\
\cmidrule(lr){2-2}\cmidrule(lr){3-5}

Metric
& Platoon
& Buildup
& Transition
& Full congestion
\\
\midrule

\multicolumn{5}{l}{\textit{$r$-Safety identification}}\\

Estimated reserve $r$
& 0.307 & 0.281 & 0.269 & 0.249 \\

Validation $R^2$ at $\hat r$
& 0.633 & 0.397 & 0.471 & 0.449 \\

Validation $R^2$ at $r=0$
& 0.616 & 0.387 & 0.464 & 0.445 \\

\midrule

\multicolumn{5}{l}{\textit{AKS representation and planning}}\\

Median reconstruction $R^2$
& 0.894 & 0.824 & 0.835 & 0.774 \\

Feasible planned references
& 92\% & 92\% & 88\% & 91\% \\

Median planning $R^2$
& 0.702 & 0.477 & 0.413 & 0.397 \\

\bottomrule
\end{tabular}
\end{table}

Together, the NGSIM results extend the two main findings of the controlled
experiments: a braking-normalized spacing reserve can be identified under
naturalistic traffic conditions, and the same AKS family can both represent
observed transitions and generate prospective references without observing
their interior trajectories. The following section places these planned
references inside a longitudinal feedback controller and evaluates their
closed-loop execution.

\section{Closed-Loop Execution of Planned Transitions}
\label{sec:execution}

The closed-loop evaluation considers two complementary settings. The controlled
platoon and NGSIM transitions are used for a same-controller ablation between
policy-only and AKS-guided execution, with the leader input, initial state,
feedback law, loop gain, and actuation limits held fixed.

The FHWA ACC/CACC experiments provide an automated-following setting with
measured follower responses under prescribed leader maneuvers. These measured
trajectories are reported alongside matched-input policy-only and AKS-guided
executions of the common feedback model.

\subsection{Execution framework and data-grounded parameterization}
\label{sec:exec_testbed}

All closed-loop experiments use the reference-tracking law introduced in
Section~\ref{sec:method}. The two execution modes differ only in how the
reference supplied to that controller is constructed.

For AKS-guided execution, the controller receives the planned finite-horizon
reference $\{s^*(t),v^*(t),a^*(t)\}$ generated by the selected AKS transition.
For policy-only tracking, the spacing policy is evaluated online at the
current follower speed,

\begin{equation}
s^*(t)=S_r\!\left(v_F(t)\right),
\qquad
v^*(t)=v_L(t),
\qquad
a^*(t)=0 .
\label{eq:policy_reference}
\end{equation}

Substituting this reference into the tracking law gives

\begin{equation}
a_{\rm cmd}^{\rm policy}(t)
=
k_s\!\left[s(t)-S_r\!\left(v_F(t)\right)\right]
+
k_v\!\left[v_L(t)-v_F(t)\right].
\label{eq:direct_policy_control}
\end{equation}

Policy-only tracking therefore regulates the instantaneous spacing and
relative-speed errors associated with the prescribed following policy. When
$r=0$, $S_r(v)$ reduces to the constant-time-gap form, and
Equation~\eqref{eq:direct_policy_control} becomes the corresponding ACC-type
spacing-feedback law. AKS-guided tracking uses the same feedback structure,
but replaces the instantaneous policy reference with a state-to-state
trajectory constructed over the prescribed horizon $[0,T]$. The two executions therefore use the same feedback controller; the distinction is whether the controller receives only the instantaneous spacing policy or an explicit finite-horizon AKS transition.

The relative weighting of spacing and speed feedback is fixed using the gain
structure adopted in previous ACC/CACC studies
\citep{milanes2014realtraffic,milanes2014modeling},

\begin{equation}
(k_s,k_v)=K(1,T_c),
\qquad
T_c=0.556~\mathrm{s}.
\label{eq:gain_split}
\end{equation}

For the controlled-platoon and NGSIM environments, $K$ is calibrated from the
corresponding recorded follower responses. For the FHWA experiments, $K$ is
calibrated using the ACC followers and the resulting value is retained for the
matched-input CACC analysis. The calibration is completed before the closed-loop ablation, and the same value of $K$ is used for both executions within each setting.

\begin{figure}[!htbp]
    \centering
    \includegraphics[width=0.75\linewidth]{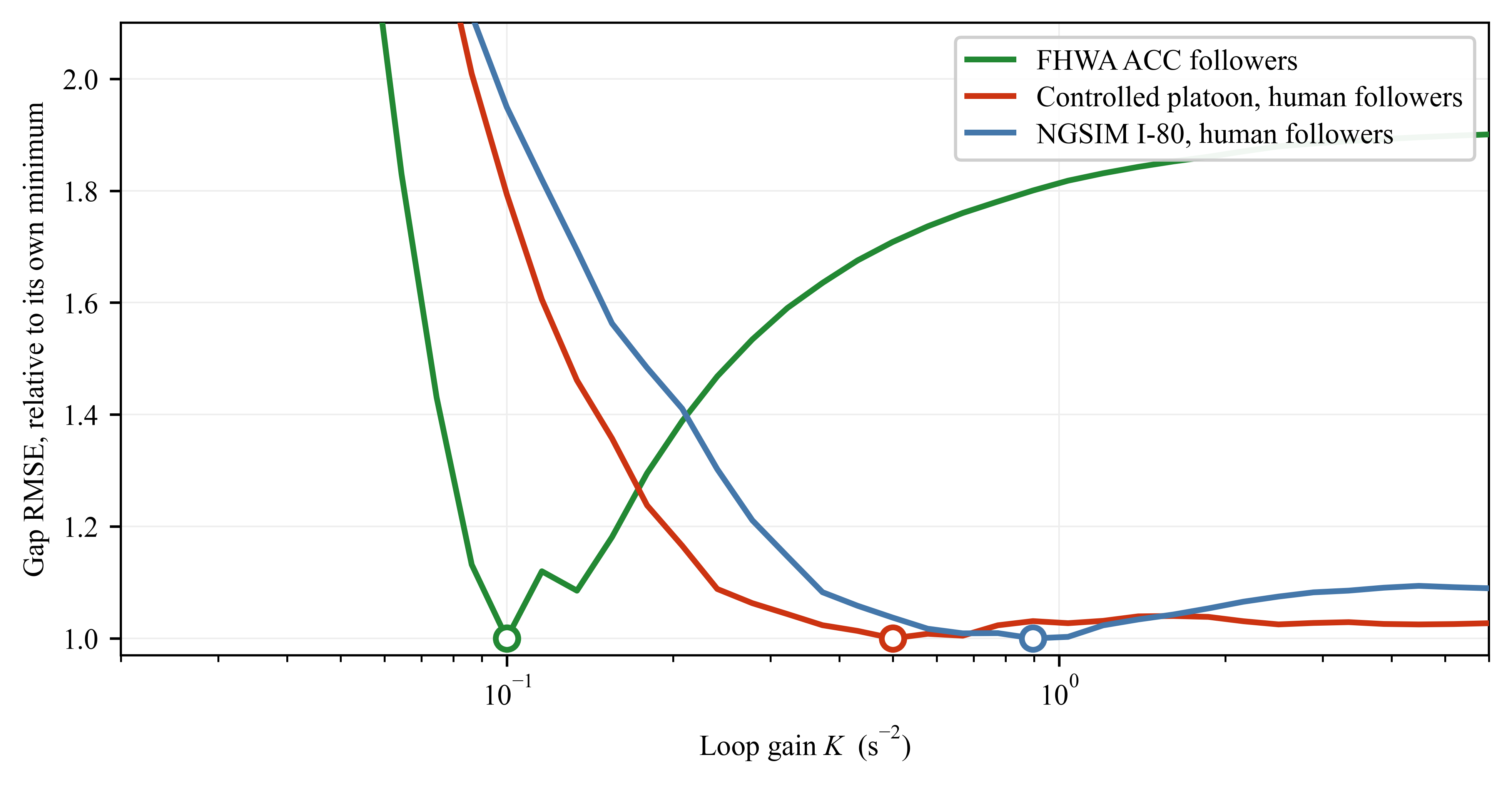}
    \caption{Calibration of the overall loop gain $K$ for the three execution
    environments. Each curve shows the gap RMSE relative to its minimum under
    the recorded leader input.}
    \label{fig:calibration}
\end{figure}

Figure~\ref{fig:calibration} shows the resulting response scales. The FHWA ACC calibration yields a substantially lower loop gain than the controlled and naturalistic human-following calibrations. These environment-specific values are used only to establish a data-grounded
execution interface; within each environment, the controller parameters remain
fixed across the policy-only and AKS-guided comparisons.

\begin{table}[!htbp]
\centering
\caption{Execution settings used in the closed-loop evaluation.}
\label{tab:controller_settings}
\begin{tabular}{llll}
\toprule
Execution setting & Spacing policy & Follower context & $K~(\mathrm{s^{-2}})$ \\
\midrule
Controlled platoon
& fitted $r$-Safety ($r=0.307$)
& controlled human trajectories
& 0.50 \\

NGSIM I-80
& period-specific $r$-Safety
& naturalistic human trajectories
& 0.90 \\

FHWA ACC
& $\tau=1.1~\mathrm{s}$, $r=0$
& ACC followers
& 0.10 \\

FHWA CACC
& $\tau=1.2~\mathrm{s}$, $r=0$
& CACC followers
& 0.10 \\
\bottomrule
\end{tabular}
\end{table}

For the controlled-platoon and NGSIM experiments, the recorded human follower
trajectory is not replayed as the closed-loop response. Only the recorded
leader motion and initial follower state are retained, and the follower is
regenerated by the common execution model under either policy-only tracking
or AKS guidance. The FHWA experiments provide a complementary automated-following setting in
which measured ACC and CACC responses can be examined alongside the
matched-input executions.

The evaluation asks whether the required spacing transition is realized within
the prescribed horizon and what motion demand is associated with that
execution. Terminal spacing error measures finite-horizon displacement closure,
while $J_a$, $J_j$, and the acceleration and jerk extrema characterize the
resulting motion. All comparisons are made under identical leader inputs,
controller parameters, and initial conditions within each event.

\subsection{Finite-horizon execution under controlled and naturalistic trajectories}
\label{sec:exec_jiang}

We first isolate the contribution of the AKS transition layer through an
internal ablation on the controlled-platoon and NGSIM
transitions examined in Sections~\ref{sec:exp_jiang} and
\ref{sec:exp_ngsim}. Feasible AKS references are available for 441 of the
477 controlled transitions and 648 of the 718 naturalistic transitions.
For each event, the recorded leader trajectory and initial follower state are
retained, while the follower motion is regenerated using the execution model
of Section~\ref{sec:exec_testbed}. Policy-only and AKS-guided
executions therefore use the same leader input, feedback law, loop gain, and
actuation limits; only the reference supplied to the controller differs. The
recorded follower response is retained in the comparison as an empirical
reference.

The comparison is made over the prescribed transition horizon $[0,T]$.
A conventional settling-time metric is not available on these records because
the next leader disturbance typically occurs only a few seconds after the
identified transition ends. The median post-transition quasi-steady interval
is approximately $5.8~\mathrm{s}$ for the controlled data and
$4.8~\mathrm{s}$ for NGSIM. The primary endpoint measure is therefore the
terminal spacing error of Eq.~\eqref{eq:tracking_errors}, evaluated at
$t=T$, which measures the residual error in realizing the displacement
required by the fitted $r$-Safety relation.

Table~\ref{tab:execution} summarizes the resulting closed-loop performance.
The first group of metrics evaluates finite-horizon transition completion,
while the second characterizes the effort and smoothness of the executed
motion. The observed follower is included only to indicate the scale of the
empirical response and is not treated as an additional control strategy.

\begin{table}[!htbp]
\centering
\caption{Internal ablation of the AKS transition layer under controlled and
naturalistic leader trajectories. The observed follower is retained as an
empirical reference. Policy-only and AKS-guided executions use the same
feedback controller, leader input, initial condition, and loop gain. Terminal
deviations quantify realization of the prescribed finite-horizon spacing
target, whereas percentages in parentheses report reductions in motion-demand
metrics relative to the policy-only execution.}
\label{tab:execution}

\small
\setlength{\tabcolsep}{2pt}
\renewcommand{\arraystretch}{1.08}

\begin{tabular}{lccc@{\hskip 8pt}ccc}
\toprule
&
\multicolumn{3}{c}{Controlled platoon ($n=441$)}
&
\multicolumn{3}{c}{NGSIM I-80 ($n=648$)}
\\
\cmidrule(lr){2-4}\cmidrule(lr){5-7}

Metric
& Observed & Policy-only & AKS-guided
& Observed & Policy-only & AKS-guided
\\
\midrule

\multicolumn{7}{l}{\textit{Finite-horizon target realization}}\\

Median terminal deviation (m)
& 3.176 & 0.668 & 0.035
& 2.078 & 0.595 & 0.045 \\

90th-percentile deviation (m)
& 8.52 & 1.92 & 0.13
& 6.39 & 1.84 & 0.39 \\

\midrule

\multicolumn{7}{l}{\textit{Motion effort and smoothness}}\\

$J_a$ (m$^2$ s$^{-3}$)
& 11.43 & 10.59 & 5.57 (\textbf{$-47.4\%$})
& 9.23 & 3.43 & 2.70 (\textbf{$-21.3\%$}) \\

$J_j$ (m$^2$ s$^{-5}$)
& 5.913 & 1.349 & 0.365 (\textbf{$-72.9\%$})
& 35.06 & 5.419 & 1.037 (\textbf{$-80.9\%$}) \\

Peak $|a|$ (m s$^{-2}$)
& 1.485 & 1.308 & 0.887 (\textbf{$-32.2\%$})
& 2.102 & 1.249 & 0.950 (\textbf{$-23.9\%$}) \\

Peak $|j|$ (m s$^{-3}$)
& 1.472 & 0.706 & 0.337 (\textbf{$-52.3\%$})
& 4.696 & 2.139 & 0.894 (\textbf{$-58.2\%$}) \\

\bottomrule
\end{tabular}
\end{table}

The transition-completion results show a consistent advantage from explicit
AKS guidance. On the controlled platoon data, the median absolute terminal gap
error decreases from $0.668~\mathrm{m}$ under policy-only tracking to
$0.035~\mathrm{m}$ with AKS guidance. The corresponding reduction on NGSIM is
from $0.595$ to $0.045~\mathrm{m}$. The separation remains in the upper tail:
the ninetieth-percentile error decreases from $1.92$ to $0.13~\mathrm{m}$ on
the controlled data and from $1.84$ to $0.39~\mathrm{m}$ in naturalistic
traffic.

The same ablation also reduces motion demand. Relative to policy-only
execution, AKS lowers $J_a$ by $47.4\%$ on the controlled data and $21.3\%$
on NGSIM, while $J_j$ decreases by $72.9\%$ and $80.9\%$, respectively.
Peak acceleration and jerk show the same overall pattern.

The two datasets therefore show the same overall closed-loop pattern under
substantially different driving conditions. Policy-only tracking responds to the
instantaneous spacing and relative-speed errors generated during the maneuver.
AKS guidance instead supplies a finite-horizon reference whose displacement has
already been matched to the prescribed terminal spacing. Under the same
feedback controller, this explicit transition guidance produces more accurate
displacement closure while requiring less acceleration effort and smoother
motion.

\subsection{Finite-horizon execution with empirical ACC/CACC trajectories}
\label{sec:exec_fhwa}

The FHWA platooning experiments provide an empirical automated-following
setting for examining finite-horizon spacing adjustment. We use transitions with ACC followers and CACC followers driven by prescribed lead-vehicle speed profiles
\citep{tiernan2017cacc}. The two configurations are analyzed separately because
their prescribed lead-vehicle dynamics differ by design. Their spacing policies
are known from the experiments rather than estimated: the commanded time gaps
are $1.1~\mathrm{s}$ for the ACC followers and $1.2~\mathrm{s}$ for the CACC
followers, corresponding to constant-time-gap spacing with $r=0$. A kinematic
consistency screen based on the relative-motion identity retains 24 ACC and
50 CACC transitions for analysis.

The ACC transitions also provide an empirical check of the common execution
model. Under policy-only tracking, the simulated follower reproduces the
measured ACC speed with a median event-level RMSE of $0.517~\mathrm{m/s}$ over
the transition horizon, and its median terminal gap error
($3.000~\mathrm{m}$) is nearly identical to that of the measured ACC follower
($2.989~\mathrm{m}$). The CACC responses are reproduced less closely, with a median speed RMSE of $1.632~\mathrm{m/s}$. We therefore retain the measured ACC and CACC
trajectories as empirical automated-following responses and report the
policy-only and AKS-guided cases as matched-input executions of the common
feedback model.

Table~\ref{tab:fhwa} summarizes the measured and simulated transitions. In
addition to terminal spacing accuracy, the table reports the fraction of the
policy-required gap change realized by the end of the speed-transition horizon.
A value of $1.00$ denotes exact spatial closure, while values below and above
one indicate under- and over-delivery, respectively.

\begin{table}[!htbp]
\centering
\caption{FHWA automated-following responses and matched-input closed-loop
executions over the prescribed speed-transition horizon. Measured ACC and CACC
responses are shown together with policy-only and AKS-guided executions generated
under the same recorded leader inputs and initial states. The gap-change fraction
is the realized gap change divided by the change required by the commanded spacing
policy; $1.00$ denotes exact closure. Entries are medians.}
\label{tab:fhwa}

\small
\setlength{\tabcolsep}{5.8pt}
\renewcommand{\arraystretch}{1.16}

\begin{tabular}{
p{2.15cm}
p{2.75cm}
ccc@{\hskip 14pt}ccc
}

\toprule
& &
\multicolumn{3}{c}{CACC followers ($n=50$)}
&
\multicolumn{3}{c}{ACC followers ($n=24$)}
\\
\cmidrule(lr){3-5}\cmidrule(lr){6-8}

Category & Metric
& Measured
& \makecell{Policy-\\only}
& \makecell{AKS-\\guided}
& Measured
& \makecell{Policy-\\only}
& \makecell{AKS-\\guided}
\\
\midrule

\multirow{2}{2.0cm}{Transition timing}
& Speed horizon $T$ (s)
& 15.3 & --- & ---
& 15.5 & --- & --- \\

& Gap settling time (s)
& 32.1 & --- & ---
& 20.4 & --- & --- \\

\midrule

\multirow{2}{2.0cm}{Spacing response}
& Gap-change fraction
& 0.12 & 1.22 & 1.00
& 1.56 & 1.22 & 1.00 \\

& Absolute terminal gap error (m)
& 7.476 & 3.416 & 0.091
& 2.989 & 3.000 & 0.046 \\

\midrule

\multirow{3}{2.0cm}{Motion demand}
& $J_a$ (m$^2$ s$^{-3}$)
& 7.042 & 6.348 & 3.956
& 5.234 & 4.816 & 5.513 \\

& $J_j$ (m$^2$ s$^{-5}$)
& 7.759 & 0.446 & 0.610
& 2.326 & 0.176 & 0.888 \\

& Peak $|a|$ (m s$^{-2}$)
& 1.312 & 1.046 & 0.921
& 1.090 & 0.820 & 1.175 \\

\bottomrule
\end{tabular}
\end{table}

The measured responses show that completing the speed transition does not
necessarily complete the associated spacing transition. By the end of the
speed horizon, the CACC followers have realized only $12\%$ of the required
gap change, whereas the ACC followers have moved beyond the required change to
a median ratio of $1.56$. The distinction is particularly clear during CACC
decelerations: the commanded spacing policy requires a median gap reduction of
$7.08~\mathrm{m}$, while the measured response instead increases the gap by
$1.00~\mathrm{m}$. Thus, close coordination of follower speed does not by
itself imply spatial closure over the same transition horizon.

The long steady portions of the FHWA experiments also allow the spacing
adjustment to be observed beyond $T$, unlike the continuously disturbed
records of Section~\ref{sec:exec_jiang}. Gap settling time is defined as the
elapsed time from transition onset to the first post-transition interval in
which the measured gap remains within $\pm2~\mathrm{m}$ of the commanded
spacing policy for at least $2~\mathrm{s}$. The CACC followers reach this
condition at a median elapsed time of $32.1~\mathrm{s}$ from transition onset,
compared with a $15.3~\mathrm{s}$ speed-transition horizon. The corresponding
times for the ACC followers are $20.4$ and $15.5~\mathrm{s}$. The measured
spacing response therefore continues beyond the speed transition in both
configurations, particularly for the CACC followers.

\begin{figure}[!htbp]
    \centering
    \includegraphics[width=\linewidth]{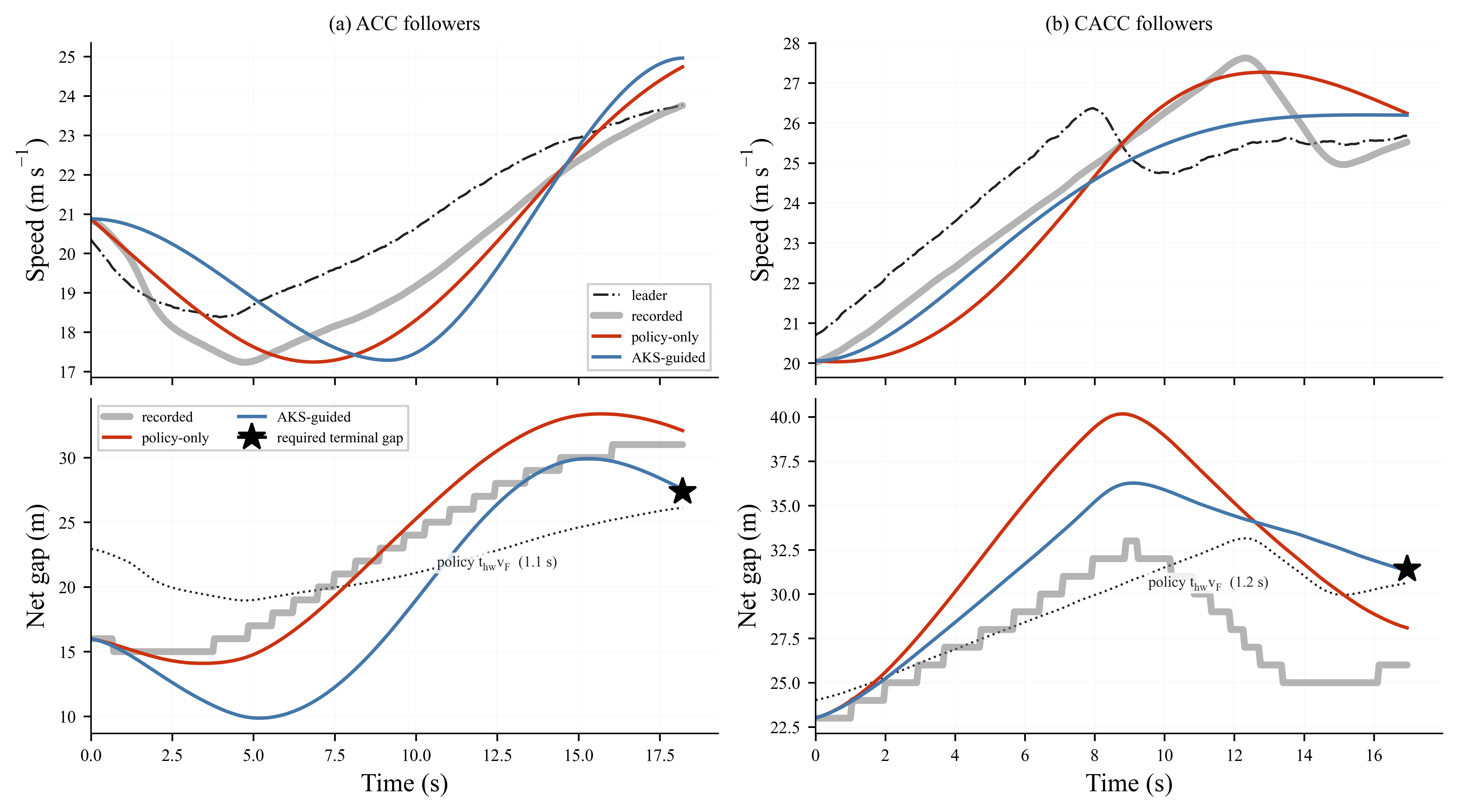}
    \caption{Representative FHWA transitions under measured automated
    following and matched-input closed-loop execution. Grey denotes the
    measured ACC or CACC follower, red policy-only tracking, and blue
    AKS-guided tracking. The marker indicates the spacing required by the
    commanded policy at the end of the speed-transition horizon.}
    \label{fig:fhwa}
\end{figure}

Under the same leader inputs and initial states, AKS guidance places the
required displacement directly into the finite-horizon reference. Within the
common execution model, the median absolute terminal gap error is reduced to
$0.091~\mathrm{m}$ for the CACC-input transitions and $0.046~\mathrm{m}$ for
the ACC-input transitions, with a gap-change fraction of $1.00$ in both cases.
Figure~\ref{fig:fhwa} illustrates this distinction for representative
transitions.

The motion-demand comparison differs across the two FHWA configurations. For
the CACC-input transitions, AKS guidance reduces $J_a$ from $6.348$ to
$3.956$ and peak acceleration from $1.046$ to $0.921~\mathrm{m/s^2}$, while
$J_j$ increases modestly from $0.446$ to $0.610$. For the ACC-input
transitions, the tighter finite-horizon closure is accompanied by increases in
$J_a$ from $4.816$ to $5.513$, $J_j$ from $0.176$ to $0.888$, and peak
acceleration from $0.820$ to $1.175~\mathrm{m/s^2}$. Thus, the effect of AKS
guidance on motion demand depends on the execution setting, whereas its effect
on finite-horizon spacing closure is consistent across both configurations.

\section{Conclusions}
\label{sec:conclusion}

This paper develops a spacing-to-transition framework for longitudinal
car-following behavior. The central idea is to connect three problems that are
usually treated separately: how spacing varies with speed, how a vehicle moves
between two speed--spacing states, and how such a transition can be used as a
reference for longitudinal control. The $r$-Safety formulation represents the
speed-dependent spacing reserve on a braking-distance scale, and the resulting
change in spacing determines the follower displacement required during a
transition. Adaptive Kinematic Smoothing (AKS) then provides a compact
analytical family for realizing that displacement over a finite horizon.

The empirical results support this construction across controlled and
naturalistic driving data. In the controlled platoon experiments, the
identified reserve is approximately $r=0.31$, while the three NGSIM periods
produce values between $0.25$ and $0.28$. The estimated reserve decreases across the three NGSIM periods as congestion develops, while the same spacing formulation remains identifiable throughout. The AKS family also captures a wide range of observed
finite transitions and can generate prospective references without using the
interior follower trajectory. Feasible references are obtained for 441 of 477
controlled transitions and 648 of 718 NGSIM transitions, with the same
low-dimensional family covering both retrospective reconstruction and
prospective planning.

The closed-loop results show why this transition layer is useful. Under the
same feedback controller, leader trajectory, and initial state, AKS-guided
tracking reduces the median terminal gap error from $0.668$ to
$0.035~\mathrm{m}$ on the controlled data and from $0.595$ to
$0.045~\mathrm{m}$ on NGSIM. In these two datasets, the improved closure is
also accompanied by lower acceleration effort and smoother motion. The FHWA
ACC/CACC experiments provide a complementary view. Their measured trajectories
show that a speed transition can finish while the corresponding spacing
adjustment is still continuing. When the same transition is supplied through
an AKS reference, the required spatial change is completed within the prescribed
horizon, although the associated motion demand depends on the execution
setting. These results place AKS as a transition-planning layer between a
spacing policy and the feedback controller rather than as a replacement for
the controller itself.

A remaining question is how broadly the $r$-Safety relation can be identified
and the AKS transition family applied. The prospective planning analysis is boundary-conditioned: the target state, transition horizon, and leader-motion information are prescribed in the present experiments, while their online determination and updating belong to a
higher-level planning layer. The present results cover controlled platoons and one naturalistic freeway dataset, but broader testing is needed across driving populations,
traffic regimes, automation systems, roadway environments, and transition
types. Establishing when a stable $r$-Safety relation can be identified, and
when the same low-dimensional AKS family remains sufficient, is an important
next step for determining the generality of the framework.

The framework also opens two extensions beyond the present study. First, the
planned transition can be evaluated jointly with vehicle energy consumption,
allowing the transition shape to be selected not only for displacement closure
and smoothness but also for fuel or energy efficiency. This provides a direct
path toward energy-aware longitudinal control that preserves the desired
spacing transition while reducing unnecessary acceleration and braking. Second, the conditional braking interpretation developed in Section~\ref{sec:conditional_safety} provides a direct interface from the calibrated spacing reserve to safety assessment. Future work can model
response delay, braking capabilities, realized spacing, and driver or controller heterogeneity to estimate event-level overlap probabilities and aggregate them across critical exposures. This would extend the present spacing and transition framework toward exposure-aware collision-risk analysis.

\appendix
\section{Extension to cross-regime spacing transitions}
\label{app:cross_regime}

The main formulation considers a transition whose initial and target states
belong to the same spacing relation. In practice, the active following regime
may also change during vehicle operation. Examples include transitions between
human driving and ACC, between independent and cooperative following, or a
handover following a change in automation or communication availability. Such
cases can be represented within the same spacing-to-displacement formulation
by allowing the endpoint states to follow different $r$-Safety relations.

Let regime $z$ be described by

\begin{equation}
S_z(v)
=
d_{0,z}
+
\tau_z v
+
r_z\frac{v^2}{2a_{\rm ref}},
\qquad
z\in\{i,j\},
\label{eq:cross_regime_policy}
\end{equation}

where the standstill clearance, time-gap parameter, and braking-normalized
reserve may differ between regimes. A transition from regime $i$ at speed
$v_1$ to regime $j$ at speed $v_2$ therefore requires

\begin{equation}
\Delta x_{F,i\rightarrow j}^{\rm req}
=
\Delta x_L
+
S_i(v_1)
-
S_j(v_2).
\label{eq:cross_regime_dx}
\end{equation}

This requirement can be decomposed as

\begin{equation}
\Delta x_{F,i\rightarrow j}^{\rm req}
=
\Delta x_L
-
\left[S_i(v_2)-S_i(v_1)\right]
+
\left[S_i(v_2)-S_j(v_2)\right].
\label{eq:cross_regime_decomp}
\end{equation}

The first term in brackets is the spacing adjustment associated with the speed
change under regime $i$. The second is the additional adjustment caused by the
regime change itself. At a common speed $v$, this cross-regime spacing
difference is

\begin{equation}
\Delta S_{i\rightarrow j}(v)
=
S_i(v)-S_j(v)
=
(d_{0,i}-d_{0,j})
+
(\tau_i-\tau_j)v
+
(r_i-r_j)\frac{v^2}{2a_{\rm ref}}.
\label{eq:cross_regime_gap}
\end{equation}

A positive $\Delta S_{i\rightarrow j}$ means that the target regime requires a
smaller spacing, so the follower must gain additional distance relative to the
leader. A negative value requires the follower to yield distance. In
particular, if the regimes differ only in the $r$-Safety reserve,

\begin{equation}
\Delta S_{i\rightarrow j}^{(r)}(v)
=
(r_i-r_j)\frac{v^2}{2a_{\rm ref}},
\label{eq:cross_regime_r_only}
\end{equation}

showing directly how a change in the identified reserve changes the required
relative displacement.

The equal-speed case illustrates the role of AKS particularly clearly. Let
$v_1=v_2=v_0$ and let the leader maintain $v_0$ over a transition of duration
$T$. Then

\begin{equation}
\Delta x_{F,i\rightarrow j}^{\rm req}
=
v_0T+\Delta S_{i\rightarrow j}(v_0).
\label{eq:cross_regime_equal_speed}
\end{equation}

Because KS1 and KS1A reduce to constant-speed motion when
$v_1=v_2$, a nonzero cross-regime spacing change requires a KS2-type temporary
speed excursion. From Eq.~\eqref{eq:ks2_displacement}, the corresponding
junction speed is

\begin{equation}
v_J
=
v_0
+
\frac{2\Delta S_{i\rightarrow j}(v_0)}{T}.
\label{eq:cross_regime_vj}
\end{equation}

A transition to a tighter regime therefore produces a temporary speed increase,
whereas a transition to a larger-spacing regime produces a temporary speed
reduction. The magnitude of the excursion decreases as a longer transition
duration is allowed.

Equations~\eqref{eq:cross_regime_dx}--\eqref{eq:cross_regime_vj} extend the
same AKS construction to planned regime changes without modifying the
analytical transition family. For a known handover, the target regime defines
the terminal spacing state and AKS constructs the corresponding finite
transition. If the regime changes unexpectedly, such as after loss of
automation or communication, the realized vehicle state at the switching time
can instead be treated as a new initial condition and the transition
recomputed under the new spacing relation. Online detection and replanning for
such unplanned changes are beyond the experiments considered in this study.

\section{Micro-to-macro implication of the $r$-Safety relation}
\label{app:micro_macro}

The $r$-Safety relation is formulated at the microscopic car-following level,
but its equilibrium spacing has a direct macroscopic interpretation. Let $L$
denote vehicle length. Because $S_r(v)$ in Eq.~\eqref{eq:r_safety} denotes the
net spacing between consecutive vehicles, the corresponding front-to-front
space headway is

\begin{equation}
H_r(v)
=
L+S_r(v)
=
L+d_0+\tau v
+r\frac{v^2}{2a_{\rm ref}} .
\label{eq:macro_headway}
\end{equation}

Under homogeneous equilibrium operation, the road space occupied by each
vehicle determines the corresponding density and flow,

\begin{equation}
k(v;r)
=
\frac{1}{H_r(v)},
\qquad
q(v;r)
=
\frac{v}{H_r(v)}
=
\frac{v}
{L+d_0+\tau v+r v^2/(2a_{\rm ref})}.
\label{eq:macro_kq}
\end{equation}

Equation~\eqref{eq:macro_kq} provides the direct microscopic-to-macroscopic
mapping of the $r$-Safety reserve. At a fixed operating speed, increasing $r$
increases the longitudinal space occupied by each vehicle and therefore
reduces the corresponding equilibrium density and flow. In particular,

\begin{equation}
\frac{\partial H_r}{\partial r}
=
\frac{v^2}{2a_{\rm ref}}
>0,
\qquad
\frac{\partial q}{\partial r}
=
-\frac{v^3}
{2a_{\rm ref}H_r(v)^2}
<0 .
\label{eq:macro_sensitivity}
\end{equation}

The effect of the reserve therefore grows with speed, consistent with its
braking-distance normalization.

The same relation can be expressed in the conventional flow--density form.
Let

\begin{equation}
A=L+d_0 .
\end{equation}

Solving Eq.~\eqref{eq:macro_headway} for the equilibrium speed associated with
density $k$ gives

\begin{equation}
V_e(k;r)
=
\begin{cases}
\dfrac{k^{-1}-A}{\tau},
& r=0, \\[9pt]
\dfrac{a_{\rm ref}}{r}
\left[
-\tau+
\sqrt{
\tau^2+
\dfrac{2r}{a_{\rm ref}}
\left(k^{-1}-A\right)
}
\right],
& r>0 .
\end{cases}
\label{eq:macro_equilibrium_speed}
\end{equation}

If the operating speed is capped by a free-flow speed $V_f$, define

\begin{equation}
k_f(r)=\frac{1}{H_r(V_f)},
\qquad
k_j=\frac{1}{A}.
\label{eq:macro_characteristic_density}
\end{equation}

The resulting equilibrium fundamental diagram can then be written as

\begin{equation}
q_r(k)
=
\begin{cases}
V_f k,
& 0\leq k\leq k_f(r), \\[4pt]
kV_e(k;r),
& k_f(r)<k\leq k_j .
\end{cases}
\label{eq:macro_fd}
\end{equation}

For $r=0$, the nonlinear reserve vanishes and
Eq.~\eqref{eq:macro_fd} reduces to the triangular relation associated with the
linear spacing rule,

\begin{equation}
q_0(k)
=
\min\left\{
V_f k,\,
\frac{1-Ak}{\tau}
\right\},
\qquad
0\leq k\leq k_j .
\label{eq:macro_triangular_fd}
\end{equation}

For $r>0$, the quadratic reserve introduces curvature into the
spacing-constrained branch. The unconstrained equilibrium flow in
Eq.~\eqref{eq:macro_kq} reaches its maximum at

\begin{equation}
v_{\rm int}^{*}
=
\sqrt{
\frac{2a_{\rm ref}(L+d_0)}{r}
},
\label{eq:macro_capacity_speed}
\end{equation}

so that, after accounting for the free-flow-speed limit,

\begin{equation}
v_c^{*}
=
\min\left\{
V_f,\,
v_{\rm int}^{*}
\right\},
\qquad
k_c^{*}
=
\frac{1}{H_r(v_c^{*})},
\qquad
q_{\max}
=
\frac{v_c^{*}}{H_r(v_c^{*})}.
\label{eq:macro_capacity}
\end{equation}

These relations give the macroscopic counterpart of the microscopic reserve
used in the main formulation. At the vehicle level, $r$ determines how much
speed-dependent spacing must be released or replenished during a finite
transition. Under homogeneous equilibrium conditions, the same reserve
determines longitudinal space consumption and therefore the density, flow, and
capacity associated with the spacing relation. The $r$-Safety formulation can
thus be carried from microscopic car-following states to their corresponding
macroscopic equilibrium traffic states.

\section*{Author contributions}

\textbf{Xuesong (Simon) Zhou:} Supervision, Conceptualization, Methodology, Writing -- review \& editing, Funding acquisition.
\textbf{Ziyi Zhang:} Conceptualization, Methodology, Software, Formal analysis, Investigation, Data curation, Writing -- original draft, Visualization.

\section*{Declaration of competing interest}

The authors declare that they have no known competing financial interests or personal relationships that could have appeared to influence the work reported in this paper.

\section*{Data availability}

The NGSIM and FHWA ACC/CACC datasets used in this study are publicly
available from their respective providers. The controlled platoon dataset
was provided by the original authors.

\section*{Acknowledgements}

The authors gratefully acknowledge Professor Rui Jiang and his research group at Beijing Jiaotong University for providing the controlled car-following platoon dataset used in this study, and Professor Zuduo Zheng of The University of Queensland for his valuable comments during the early development of this
work. This work was supported by the U.S. National Science Foundation (NSF) under Grant No.~TIP-2303748, ``POSE: Phase II: CONNECT: Consortium of Open-source Planning Models for Next-generation Equitable and Efficient Communities and Transportation.''

\section*{Declaration of generative AI and AI-assisted technologies in the\\ manuscript preparation process}

During the preparation of this work, the authors used ChatGPT and Claude
to assist with language refinement and manuscript organization. After using
these tools, the authors reviewed and edited the content as needed and take
full responsibility for the content of the publication.

\bibliographystyle{plainnat}
\bibliography{references}

\end{document}